\documentclass[a4paper]{amsart}

\usepackage{amsmath, amssymb, amsfonts, amsthm, amscd, amsopn, amstext, amsxtra, euscript}

\usepackage{graphicx} 
\usepackage{xcolor}   

\usepackage[utf8]{inputenc}

\usepackage{enumitem}
\usepackage{booktabs}

\usepackage[all]{xy}

\usepackage[citation-order, msc-links]{amsrefs}

\usepackage{hyperref}
\hypersetup{
  colorlinks   = true,
  urlcolor     = blue,
  linkcolor    = blue,
  citecolor    = blue
}

\usepackage[capitalise]{cleveref}

\newtheorem{thm}{Theorem}[section]
\newtheorem{lem}[thm]{Lemma}
\newtheorem{prop}[thm]{Proposition}
\newtheorem{cor}[thm]{Corollary}

\theoremstyle{definition}
\newtheorem{defn}{Definition}[section]
\newtheorem{exa}[defn]{Example}

\theoremstyle{remark}
\newtheorem{rem}[defn]{Remark}

\crefname{thm}{Thm.}{Thms.}
\Crefname{thm}{Thm.}{Thms.}
\crefname{prop}{Prop.}{Props.}
\Crefname{prop}{Prop.}{Props.}
\crefname{lem}{Lem.}{Lems.}
\Crefname{lem}{Lem.}{Lems.}
\crefname{cor}{Cor.}{Cors.}
\Crefname{cor}{Cor.}{Cors.}
\crefname{rem}{Rem.}{Rems.}
\Crefname{rem}{Rem.}{Rems.}
\crefname{defn}{Def.}{Defs.}

\newcommand{\cP}{\mathcal{P}}

\newcommand{\F}{\mathbb{F}}

\newcommand{\bP}{\mathbb{P}}

\newcommand{\w}{\mathbf{w}}

\newcommand{\Aut}{\operatorname{Aut}}

\DeclareMathOperator\Supp{Supp}

\DeclareMathOperator\ord{ord}

\DeclareMathOperator\sing{Sing}

\DeclareMathOperator\Proj{Proj}

\DeclareMathOperator\Fix{Fix}

\newcommand{\WPw}{\mathbb{WP}_{\mathbf{w}}^{n}}

\newcommand{\I}{\mathbf{1}}

\usepackage{lipsum} 

\usepackage{amsaddr}

\providecommand{\F}{\mathbb{F}}
\providecommand{\bP}{\mathbb{P}}
\providecommand{\cP}{\mathcal{P}}
\providecommand{\w}{\mathbf{w}}
\providecommand{\Proj}{\operatorname{Proj}}
\providecommand{\Aut}{\operatorname{Aut}}
\providecommand{\Supp}{\operatorname{Supp}}
\providecommand{\sing}{\operatorname{Sing}}
\providecommand{\WPw}{W\bP^n_\w}
\providecommand{\ord}{\operatorname{ord}}
\providecommand{\Fix}{\operatorname{Fix}}
\providecommand{\Br}{\operatorname{Br}}
\providecommand{\im}{\operatorname{im}}
\providecommand{\I}{\mathbf{1}}
\providecommand{\Gm}{\mathbb{G}_m}
\providecommand{\bbA}{\mathbb{A}}
\providecommand{\Frob}{\operatorname{Frob}}
\providecommand{\sm}{\operatorname{sm}}
\providecommand{\Mass}{\operatorname{Mass}}
\providecommand{\isoclasses}{\pi_0}
\providecommand{\ZZ}{\mathbb{Z}}
\providecommand{\bbm}{\boldsymbol{\mu}}
\providecommand{\Ztw}{Z_{\mathrm{tw}}}

\title[Twists of weighted projective stacks]{Orbit counts and twist zeta functions of weighted projective stacks over finite fields}

\author{S. Salami}
\address{Institute of Mathematics and Statistics, \\
Rio de Janeiro State University, \\
Maracan\~a, Rio de Janeiro, 20950-000, RJ, Brazil}
\email{Sajad.salami@ime.uerj.br}

\author{T. Shaska}
\address{Department of Mathematics and Statistics, \\
Oakland University, \\
Rochester, MI 48309, USA}
\email{shaska@oakland.edu}


\begin{document}

\begin{abstract}
Let $\F_q$ be a finite field and $\w=(w_0,\dots,w_n)$ a vector of positive integer weights.  Several finite-field counts attached to the weighted projective space $\bP^n_\w$ are easily conflated: the coarse rational-point count and the stacky mass are both weight-independent, whereas the number $A_\w(q)$ of $\F_q^\times$-orbits on nonzero $\F_q$-representatives for the weighted action depends on the weights.  We prove the closed formula
\[
        A_\w(q)=\sum_{\emptyset\ne S\subseteq\{0,\dots,n\}}(q-1)^{|S|-1}\gcd(k_S,q-1),
        \qquad k_S=\gcd\{w_i:i\in S\},
\]
and identify $A_\w(q)$ intrinsically as the number of $\F_q$-isomorphism classes of the weighted projective stack $\cP_\w=[(\bbA^{n+1}\setminus\{0\})/\Gm]$---equivalently, the number of $\F_q$-twists lying over the coarse points---the discrepancy from the coarse count being governed by the Kummer groups $\F_q^\times/(\F_q^\times)^{k_S}$.  We read the behaviour of $A_\w$ under reduction of the weight vector through the $\F_q$-cohomology of the associated $\bbm_d$-gerbe, and prove that the twist zeta function $\Ztw(\cP_\w,t)=\exp\bigl(\sum_{r\ge1}A_\w(q^r)t^r/r\bigr)$ is rational with multiplicity spectrum independent of $q$, admitting a single global functional equation precisely when the weights share a common prime-to-$p$ part---for reduced $\w$, precisely when the twist theory is trivial.  For weighted diagonal hypersurfaces and same-degree pairs in the split regime we compute the twist zeta function of the substack explicitly, with reciprocal roots given by Gauss sums and, for intersections, by Frobenius eigenvalues of superelliptic curves.
\end{abstract}

\maketitle

\section{Introduction}

Weighted projective spaces are basic singular toric varieties; they arise as ambient spaces for weighted hypersurfaces, in the construction of stacks with quotient singularities, and in coding theory \cite{Mori1975,Hosgood,aubry0,AubryPerret2025,shaska2025}.  For a field $k$ and positive integers $\w=(w_0,\dots,w_n)$, the weighted projective space $\bP^n_\w=\Proj k[x_0,\dots,x_n]$ with $\deg x_i=w_i$ carries an associated quotient stack $\cP_\w$ whose coarse space is $\bP^n_\w$.  The stack records the isotropy of each point---the cyclic stabilizer attached to its support---which the coarse space forgets.

Over a finite field $\F_q$ several counts are attached to $\bP^n_\w$, and they are easily conflated.  The coarse rational-point count is the weight-independent value $\#\bP^n_\w(\F_q)=(q^{n+1}-1)/(q-1)$, a standard toric computation; this value is classical, is unchanged under reduction $\w\mapsto\w/d$, and is not revised here, as established by Aubry--Perret \cite{AubryPerret2025}, Phillips \cite{Phillips2024}, and Nardi--San-Jos\'e \cite{NardiSanJose}.  The stacky mass $\sum_{[x]}\#\Aut(x)^{-1}$ turns out to equal this same value.  By contrast, the orbit count $A_\w(q)$---the number of $\F_q^\times$-orbits of nonzero $\F_q$-vectors under the weighted action---is a genuinely weight-dependent invariant, \emph{distinct from the coarse point count}.  The purpose of this note is to isolate $A_\w(q)$, compute it, and explain its meaning.

The discrepancy between $A_\w(q)$ and the coarse count is a twist phenomenon.  A single coarse $\F_q$-point may carry several inequivalent $\F_q$-representatives under the smaller group $\F_q^\times$, which become equivalent only after extension of scalars; the local obstruction over a point of support $S$ is a Kummer group of order $\gcd(k_S,q-1)$ (Section~\ref{sec:orbit-count}).  Summing over coarse points yields the closed formula
\[
        A_\w(q)=\sum_{\emptyset\ne S\subseteq T}(q-1)^{|S|-1}\gcd(k_S,q-1),
        \qquad k_S=\gcd\{w_i:i\in S\},
\]
and identifies $A_\w(q)$ with the number of $\F_q$-isomorphism classes of $\cP_\w$, that is, the total count of $\F_q$-twists over the coarse points.

Our contribution is to

(i) separate cleanly the coarse point count and the stacky mass---both weight-independent---from the orbit count $A_\w$, which is weight-dependent, together with the explicit Kummer mechanism producing the discrepancy;

(ii) give the intrinsic interpretation of $A_\w(q)$ as the isomorphism-class count of the groupoid $\cP_\w(\F_q)$, assembled from the local twist groups, and read its behaviour under reduction of the weight vector, $\w\mapsto\w/d$ with $d=\gcd(w_0,\dots,w_n)$---an operation that preserves the coarse variety but changes the stack by a $\bbm_d$-gerbe---through the $\F_q$-cohomology of this gerbe;

(iii) introduce the twist zeta function $\Ztw(\cP_\w,t)$, prove its rationality through an explicit cyclotomic product, and settle its symmetry theory: a sector decomposition of the orbit count yields the multiplicity spectrum in closed form, revealing it to be nonincreasing and independent of $q$, and a single global functional equation holds precisely when the weights share a common prime-to-$p$ part---for reduced $\w$, precisely when the twist theory is trivial; and

(iv) extend the orbit count to weighted complete intersections and compute the twist zeta function of diagonal hypersurfaces and of same-degree pairs in the split regime, where it is rational with reciprocal roots given by Gauss sums and, for intersections, by Frobenius eigenvalues of superelliptic curves.

We regard $A_\w$ as the \emph{local}, fixed-field counterpart of the \emph{global} counts of points of bounded height on weighted projective stacks studied by Phillips \cite{Phillips2024} and, in the general framework of heights on stacks, by Ellenberg, Satriano, and Zureick-Brown \cite{ESZB}; the same $\bbm$-torsor mechanism produces, over number fields, the arithmetic sparsity obstruction of \cite{shaska-sparsity}.  Heights and entropy over finite fields are taken up in the companion paper \cite{ss-heights}.

The paper is organized as follows.  Section~\ref{sec:rationality} treats rational points on the coarse space; Section~\ref{sec:orbit-count} proves the orbit-count formula and its twist interpretation; Section~\ref{sec:normalization} studies reduction of weights through the gerbe cohomology; Section~\ref{sec:strata} gives the smooth/singular decomposition; Section~\ref{sec:zeta} proves rationality of the twist zeta function; Section~\ref{sec:fe} computes the multiplicity spectrum and classifies the weight vectors admitting a global functional equation; Section~\ref{sec:hypersurfaces} treats weighted complete intersections and the twist zeta functions of diagonal hypersurfaces; Section~\ref{sec:conclusion} gives concluding remarks, and Appendix~\ref{app:tables} records numerical tables.

\medskip

\noindent\textbf{Notation.}  Throughout this paper, $\F_q$ is a finite field with $q=p^a$ elements and $\overline{\F}_q$ an algebraic closure. We always denote by  $T=\{0,\dots,n\}$.

\section{Rational points on weighted projective spaces}\label{sec:rationality}

 The weighted action of $\lambda\in\overline{\F}_q^\times$ on a vector $x=(x_0,\dots,x_n)$ is
\[
        \lambda\cdot x=(\lambda^{w_0}x_0,\dots,\lambda^{w_n}x_n),
\]
and the support of $x$ is $\Supp(x)=\{i:x_i\ne0\}$.  
For  $T=\{0,\dots,n\}$   and  $\emptyset\ne S\subseteq T$ we set $k_S=\gcd\{w_i:i\in S\}$, and $d=\gcd(w_0,\dots,w_n)=k_T$.  The singular locus of $\bP^n_\w$ is $\sing(\bP^n_\w)$ and the smooth (weak) locus is $\WPw$; $\ord_m(q)$ denotes the multiplicative order of $q$ modulo $m$.

We begin by fixing terminology with care, since three a priori distinct counts are involved and conflating them is a known source of error.  The multiplier relating two representatives of the same point of $\bP^n_\w$ need not lie in $\F_q$, even when both representatives have coordinates in $\F_q$; this single fact separates the counts below.

\begin{defn}\label{def:coarse-point}
Let $x\sim_{\F_q}x'$ mean $x'=\lambda\cdot x$ for some $\lambda\in\F_q^\times$ (weighted action), and $x\sim_{\overline{\F}_q}x'$ mean $x'=\lambda\cdot x$ for some $\lambda\in\overline{\F}_q^\times$, where the multiplier $\lambda$ is allowed to lie outside $\F_q$.  We distinguish:

\begin{enumerate}[label=\textup{(\alph*)}]
\item the \emph{coarse points} of $\bP^n_\w$ over $\F_q$, namely the $\F_q$-points of the coarse scheme 
\[
\bP^n_\w(\F_q)=\bigl(\Proj\F_q[x_0,\dots,x_n]\bigr)(\F_q),
\]
 equivalently the $\Frob_q$-fixed geometric points; a geometric point is \emph{Frobenius-invariant} if $\Frob_q(P)=P$;

\item whether such a point \emph{admits an $\F_q$-representative}, i.e.\ is represented by some $x\in\F_q^{n+1}\setminus\{0\}$;

\item the \emph{orbit count} 
\[
A_\w(q)=\#\bigl((\F_q^{n+1}\setminus\{0\})/\!\sim_{\F_q}\bigr),
\]
 the number of $\F_q$-vectors up to multipliers in $\F_q^\times$ \emph{only} (Definition~\ref{def:Aw}).
\end{enumerate}

A coarse point is a $\sim_{\overline{\F}_q}$-class; the orbit count refines it by the smaller relation $\sim_{\F_q}$.  In particular $A_\w(q)$ is \emph{not} the coarse point count $\#\bP^n_\w(\F_q)$: a single coarse point may split into several $\sim_{\F_q}$-classes, identified only by a multiplier $\lambda\in\overline{\F}_q^\times\setminus\F_q^\times$.
\end{defn}

\begin{lem}\label{lem:rationality}
For a geometric point $P\in\bP^n_\w(\overline{\F}_q)$ the following are equivalent:
\begin{enumerate}[label=\textup{(\roman*)}]
\item $P$ is Frobenius-invariant;
\item $P$ is an $\F_q$-point of the coarse variety;
\item $P$ admits a representative in $\F_q^{n+1}\setminus\{0\}$.
\end{enumerate}
\end{lem}

\begin{proof}
The equivalence of (i) and (ii) is the usual Galois criterion for rational points of a variety over a finite field.  It remains to prove that a rational point of the coarse weighted projective space has an $\F_q$-representative.
Fix a nonempty support $S$.  The corresponding locally closed toric stratum is the quotient of the split torus $\Gm^S$ by the one-dimensional subtorus
\[
        \Gm\longrightarrow \Gm^S,
        \qquad
        \lambda\longmapsto (\lambda^{w_i/k_S})_{i\in S}.
\]
The exponents $w_i/k_S$ have greatest common divisor $1$, so this is a closed embedding of split tori.  Hence there is an exact sequence of split tori over $\F_q$
\[
        1\longrightarrow \Gm\longrightarrow \Gm^S\longrightarrow \mathbb T_S\longrightarrow 1,
\]
where $\mathbb T_S$ is the torus underlying the stratum.  Taking Galois cohomology and using Hilbert 90, equivalently Lang's theorem for $\Gm$, gives $H^1(\F_q,\Gm)=0$.  Therefore 
\[
\Gm^S(\F_q)\to\mathbb T_S(\F_q)
\]
 is surjective, and every rational point in the stratum lifts to an $S$-supported vector with coordinates in $\F_q^\times$.  This proves (ii)$\Rightarrow$(iii).  The implication (iii)$\Rightarrow$(i) is immediate, since Frobenius fixes every coordinate of such a representative.
\end{proof}

The lifting in (ii)$\Rightarrow$(iii) is special to finite fields: it rests on $H^1(\F_q,\Gm)=0$.  Over a number field the analogous map fails to be surjective on rational points, and the resulting $\bbm$-torsor obstruction thins the count of points of bounded weighted height \cite{shaska-sparsity}; over $\F_q$ the same torsors obstruct nothing but instead multiply rational representatives, and that multiplication is precisely the invariant $A_\w(q)$ of this paper.

\begin{cor}\label{cor:coarse-count}
For every vector of positive weights $\w$ one has
\[
        \#\bP^n_\w(\F_q)=\sum_{\emptyset\ne S\subseteq T}(q-1)^{|S|-1}=\frac{q^{n+1}-1}{q-1}.
\]
\end{cor}

\begin{proof}
By the proof of \cref{lem:rationality}, the stratum of support $S$ is the split torus $\mathbb T_S\simeq\Gm^{|S|-1}$, with $(q-1)^{|S|-1}$ rational points.  Summing over all nonempty supports gives 
\[\sum_{m=1}^{n+1}\binom{n+1}{m}(q-1)^{m-1}=(q^{n+1}-1)/(q-1).\]
\end{proof}

\begin{rem}\label{rem:coarse-known}
The coarse count $\#\bP^n_\w(\F_q)=(q^{n+1}-1)/(q-1)$ is classical, weight-independent, and unchanged under reduction $\w\mapsto\w/d$; see Aubry--Perret \cite{AubryPerret2025}, Phillips \cite{Phillips2024}, and the appendix of Nardi--San-Jos\'e \cite{NardiSanJose}, where the weight-independence and the equality of $\F_q$-point sets $\bP^n_\w(\F_q)=\bP^n_{\w/d}(\F_q)$ are established.  We do not revise this count.  The invariant studied in this paper is the orbit count $A_\w(q)$ of \cref{def:coarse-point}(c), which is a different quantity and is the one sensitive to the weights; the coincidence of the coarse count with the stacky mass is recorded in \cref{prop:mass-equals-coarse}.
\end{rem}

\begin{exa}\label{exa:strict-representative}
Let $q=3$ and $\w=(1,2)$, and choose $i\in\F_9$ with $i^2=-1$.  The geometric point $[i,1]$ is Frobenius-invariant, since 
\[
(i^3,1^3)=(-i,1)=(-1)\cdot(i,1)
\]
 for the weighted action; the multiplier $-1=(-i)^2$ is realized only by the non-rational scalar $-i\in\F_9\setminus\F_3$.  Rescaling by $-i$ gives $[i,1]=[1,-1]$, an $\F_3$-representative.  
 
 Thus $[i,1]$ is a single coarse $\F_3$-point: a non-rational representative must not be mistaken for the absence of a rational one, nor counted as an additional coarse point.  The point $[0:y]$ with $y\in\F_q^\times$ is likewise the single coarse point $[0:1]$, since $y=\lambda^2$ for some $\lambda\in\overline{\F}_q^\times$ rescales $(0,y)$ to $(0,1)$.
\end{exa}

\section{Orbit counts and twist groups}\label{sec:orbit-count}

The coarse count of \cref{cor:coarse-count} uses geometric equivalence under $\overline{\F}_q^\times$.  We now restrict the scaling parameter to $\F_q^\times$.

\begin{defn}\label{def:Aw}
The finite-field orbit count attached to $\w$ is
\[
\begin{split}
        A_\w(q) & =\#\left((\F_q^{n+1}\setminus\{0\})/\F_q^\times\right), \; \text{ where } \\
        \lambda\cdot(x_0,\dots,x_n) & =(\lambda^{w_0}x_0,\dots,\lambda^{w_n}x_n),
\end{split}        
\]
for  $\lambda\in\F_q^\times$.
\end{defn}

Then we have the following well known result. 

\begin{lem}[Burnside formula]\label{lem:burnside}
For $\lambda\in\F_q^\times$ set $N(\lambda)=\#\{i\in T:\lambda^{w_i}=1\}$.  Then
\[
        A_\w(q)=\frac1{q-1}\sum_{\lambda\in\F_q^\times}\bigl(q^{N(\lambda)}-1\bigr).
\]
\end{lem}

\begin{proof}
Burnside's lemma applied to the $\F_q^\times$-action gives 
\[
A_\w(q)=\tfrac1{q-1}\sum_{\lambda}\#\Fix(\lambda)
\]
 on $\F_q^{n+1}\setminus\{0\}$.  A coordinate $x_i$ is fixed by $\lambda$ if and only if $x_i=0$ or $\lambda^{w_i}=1$, so the fixed locus in $\F_q^{n+1}$ has cardinality $q^{N(\lambda)}$, and deleting the zero vector gives $q^{N(\lambda)}-1$.
\end{proof}

\begin{lem}[Support formula]\label{lem:orbit-formula}
One has
\[
        A_\w(q)=\sum_{\emptyset\ne S\subseteq T}(q-1)^{|S|-1}\gcd(k_S,q-1).
\]
\end{lem}

\begin{proof}
Stratify by support.  The number of vectors of support exactly $S$ is $(q-1)^{|S|}$; each has stabilizer $\{\lambda\in\F_q^\times:\lambda^{w_i}=1\ \forall i\in S\}=\bbm_{k_S}(\F_q)$, of order $\gcd(k_S,q-1)$, hence orbit size $(q-1)/\gcd(k_S,q-1)$.  The stratum therefore contributes $(q-1)^{|S|-1}\gcd(k_S,q-1)$ orbits, and summing over nonempty $S$ gives the formula.
\end{proof}

\begin{prop}\label{prop:burnside-support-agree}
The Burnside formula of \cref{lem:burnside} and the support formula of \cref{lem:orbit-formula} agree.
\end{prop}

\begin{proof}
For each $\lambda\in\F_q^\times$,
\[
        q^{N(\lambda)}
        =\prod_{i=0}^n\bigl(1+(q-1)\,\I(\lambda^{w_i}=1)\bigr)
        =\sum_{S\subseteq T}(q-1)^{|S|}\,\I(\lambda^{k_S}=1),
\]
with the convention that the term $S=\emptyset$ equals $1$.  Hence
\[
\begin{aligned}
        A_\w(q)
        &=\frac1{q-1}\sum_{\lambda\in\F_q^\times}\bigl(q^{N(\lambda)}-1\bigr)\\
        &=\frac1{q-1}\sum_{\emptyset\ne S\subseteq T}(q-1)^{|S|}\,\#\{\lambda\in\F_q^\times:\lambda^{k_S}=1\}\\
        &=\sum_{\emptyset\ne S\subseteq T}(q-1)^{|S|-1}\gcd(k_S,q-1),
\end{aligned}
\]
using $\#\{\lambda\in\F_q^\times:\lambda^{k_S}=1\}=\gcd(k_S,q-1)$.
\end{proof}

We record the joint statement for reference.

\begin{thm}[Orbit count formula]\label{thm:orbit-count}
For $\lambda\in\F_q^\times$ set $N(\lambda)=\#\{i\in T:\lambda^{w_i}=1\}$.  Then
\[
        A_\w(q)
        =\frac1{q-1}\sum_{\lambda\in\F_q^\times}\bigl(q^{N(\lambda)}-1\bigr)
        =\sum_{\emptyset\ne S\subseteq T}(q-1)^{|S|-1}\gcd(k_S,q-1).
\]
\end{thm}

\begin{proof}
Combine \cref{lem:burnside,lem:orbit-formula,prop:burnside-support-agree}.
\end{proof}

\begin{cor}\label{cor:extensions}
For every $r\ge1$, one has
\[
        A_\w(q^r)=\sum_{\emptyset\ne S\subseteq T}(q^r-1)^{|S|-1}\gcd(k_S,q^r-1).
\]
\end{cor}


\begin{exa}\label{exa:P12}
For $\bP_{(1,2)}$ one obtains
\[
        A_{(1,2)}(q)=
        \begin{cases}
        q+2,& q\text{ odd},\\
        q+1,& q\text{ even}.
        \end{cases}
\]
Indeed, the support $\{0,1\}$ contributes $q-1$, the support $\{0\}$ contributes $1$, and the support $\{1\}$ contributes $\gcd(2,q-1)$.  The coarse count is always $q+1$.
\end{exa}

\begin{exa}[Ordinary projective space]\label{exa:ordinary-projective}
If $\w=(1,\dots,1)$ then $k_S=1$ for every nonempty $S$, so 
\[
A_{(1,\dots,1)}(q)=\sum_{\emptyset\ne S\subseteq T}(q-1)^{|S|-1}=1+q+\cdots+q^n.
\]
  Thus the orbit count and the coarse point count coincide for ordinary projective space.
\end{exa}

\begin{exa}[The weights $(1,2,3,5)$]\label{exa:orbit-1235}
For $\w=(1,2,3,5)$ the only nontrivial one-coordinate stabilizers are $\bbm_2,\bbm_3,\bbm_5$, while all supports of cardinality $\ge2$ have gcd $1$.  Hence
\[
        A_{(1,2,3,5)}(q)= q^3+q^2+q-2+\gcd(2,q-1)+\gcd(3,q-1)+\gcd(5,q-1).
\]
\end{exa}

\subsection{Twist groups and stacky interpretation}\label{sec:twists}

We now explain the intrinsic meaning of the factor $\gcd(k_S,q-1)$.

\begin{lem}\label{lem:kummer}
There is a canonical isomorphism of finite abelian groups
\[
        H^1(\F_q,\bbm_{k_S})\simeq \F_q^\times/(\F_q^\times)^{k_S},
\]
and this group has order $\#H^1(\F_q,\bbm_{k_S})=\gcd(k_S,q-1)$.
\end{lem}

\begin{proof}
The Kummer sequence
\[
1\to\bbm_{k_S}\to\Gm\xrightarrow{(\cdot)^{k_S}}\Gm\to1
\]
is exact for the fppf (faithfully flat, finite presentation) topology---which is needed rather than the \'etale topology when $p\mid k_S$, since then $\bbm_{k_S}$ is non-reduced---and together with $H^1(\F_q,\Gm)=0$ (Hilbert~90) it gives
\[
H^1(\F_q,\bbm_{k_S})\simeq\F_q^\times/(\F_q^\times)^{k_S}.
\]
   Since $\F_q^\times$ is cyclic of order $q-1$, this quotient has order $\gcd(k_S,q-1)$.
\end{proof}

\begin{prop}\label{prop:twists-over-point}
Let $P\in\bP^n_\w(\F_q)$ have support $S$.  The set of $\F_q^\times$-orbits of $\F_q$-representatives lying over the coarse point $P$ is naturally a torsor under $H^1(\F_q,\bbm_{k_S})$.  Consequently the number of such orbits is $\gcd(k_S,q-1)$.
\end{prop}

\begin{proof}
Choose an $\F_q$-representative $v$ of $P$, which exists by \cref{lem:rationality}.  Over $\overline{\F}_q$ all representatives of the same geometric point have the form $\alpha\cdot v$, $\alpha\in\overline{\F}_q^\times$, and the stabilizer of $v$ is $\bbm_{k_S}$.  Rational forms of $v$ modulo rational equivalence are therefore classified by $H^1(\F_q,\bbm_{k_S})$; as $\bbm_{k_S}$ is abelian this set is a group, of cardinality $\gcd(k_S,q-1)$ by \cref{lem:kummer}.
\end{proof}

\begin{thm}\label{thm:stack-classes}
Let $\cP_\w=[(\bbA^{n+1}\setminus\{0\})/\Gm]$.  Then
\[
        A_\w(q)=\#\isoclasses(\cP_\w(\F_q))
        =\sum_{P\in\bP^n_\w(\F_q)}\#H^1(\F_q,\Aut_{\cP_\w}(P)).
\]
If $P$ has support $S$, then $\Aut_{\cP_\w}(P)=\bbm_{k_S}$ and the summand is $\gcd(k_S,q-1)$.
\end{thm}

\begin{proof}
An object of $\cP_\w(\F_q)$ is a $\Gm$-torsor over $\F_q$ with an equivariant morphism to $\bbA^{n+1}\setminus\{0\}$.  Since $H^1(\F_q,\Gm)=0$ every such torsor is trivial, so isomorphism classes of objects are exactly $\F_q^\times$-orbits on $\F_q^{n+1}\setminus\{0\}$, giving the first equality.  Grouping objects by their image in the coarse space and applying \cref{prop:twists-over-point} over each point of support $S$ gives the second.
\end{proof}

\begin{rem}\label{rem:not-mass}
The count in \cref{thm:stack-classes} is an \emph{unweighted} isomorphism-class count, not the stacky mass 
\[
\Mass(\cP_\w(\F_q))=\sum_{[x]}\#\Aut(x)(\F_q)^{-1}.
\]
  We study the former because it is precisely the orbit count of rational representatives, the count that detects Kummer twists over coarse points.
\end{rem}

\begin{prop}\label{prop:mass-equals-coarse}
The stacky mass of $\cP_\w(\F_q)$ equals the coarse point count:
\[
        \Mass(\cP_\w(\F_q))
        =\sum_{[x]\in\isoclasses(\cP_\w(\F_q))}\frac1{\#\Aut(x)(\F_q)}
        =\#\bP^n_\w(\F_q)
        =\frac{q^{n+1}-1}{q-1}.
\]
In particular the mass is weight-independent, and among the coarse count, the mass, and the orbit count, only $A_\w(q)$ is sensitive to the weights.
\end{prop}

\begin{proof}
Group the isomorphism classes by their image in the coarse space.  Over a coarse point of support $S$ there are, by \cref{prop:twists-over-point}, exactly $\gcd(k_S,q-1)$ classes, each with automorphism group $\bbm_{k_S}(\F_q)$ of order $\gcd(k_S,q-1)$; their total mass is therefore $1$.  Summing over $\bP^n_\w(\F_q)$ and using \cref{cor:coarse-count} gives the claim.
\end{proof}

\section{Behaviour under reduction of the weight vector}\label{sec:normalization}

Let $d=\gcd(w_0,\dots,w_n)$ and $\w'=\w/d$.  We call $\w'$ the \emph{reduced} weight vector and say $\w$ is reduced if $d=1$.  Then $\bP^n_\w\cong\bP^n_{\w'}$ as coarse varieties, but the stacks differ: $\cP_\w$ is a generic $\bbm_d$-gerbe over $\cP_{\w'}$.

\begin{rem}
The reduction $\w\mapsto\w/d$ studied here is an operation on the \emph{weight vector} and should not be confused with the normalization of a \emph{point} to a representative of weighted greatest common divisor $1$ in the sense of \cite{BeshajGutierrezShaska,s-sh}.  Reduction changes the weight vector and the stack---by a $\bbm_d$-gerbe---while preserving the coarse variety; point normalization selects a distinguished representative of a fixed point and leaves the weights unchanged.  We reserve ``reduction'' for the former throughout, and ``well-formed'' for the stronger condition $\gcd(w_0,\dots,\widehat{w_i},\dots,w_n)=1$ for all $i$, used in \cref{sec:strata}.
\end{rem}

\begin{prop}\label{prop:mu-d-gerbe}
	Let $k$ be a field, $\w=(w_0,\dots,w_n)$, $d=\gcd(w_0,\dots,w_n)$, and $\w'=\w/d$.
	Then there is a natural morphism
	\(
	\pi:\mathcal P_\w\longrightarrow \mathcal P_{\w'}
	\)
	which is a $\bbm_d$-gerbe.  Equivalently, after base change to an algebraic closure of $k$, every geometric fibre of $\pi$ is equivalent to $B\bbm_d$.  In particular the induced map on coarse moduli spaces is an isomorphism $\mathbb P_\w^n \cong \mathbb P_{\w'}^n$.  Moreover, if $x\in\mathcal P_\w(\overline{k})$ has support $S$, its stabilizer in $\mathcal P_\w$ is $\bbm_{k_S}$, whereas the stabilizer of its image in $\mathcal P_{\w'}$ is $\bbm_{k_S/d}$.  Thus reduction by $d$ rigidifies the generic subgroup $\bbm_d$ of the stabilizers.
\end{prop}

\begin{proof}
	Let $U=\mathbb A_k^{n+1}\setminus\{0\}$.  The stack $\mathcal P_\w$ is the quotient stack for the action $\lambda\cdot(x_0,\dots,x_n)=(\lambda^{w_0}x_0,\dots,\lambda^{w_n}x_n)$, $\lambda\in\Gm$.  Since $w_i=d w_i'$ for every $i$, this action factors through $[d]:\Gm\to\Gm$, $\lambda\mapsto\lambda^d$, because $\lambda^{w_i}=(\lambda^d)^{w_i'}$.  The kernel of $[d]$ is $\bbm_d$, so we have an exact sequence of group schemes
	\[
	1\longrightarrow \bbm_d\longrightarrow \Gm\xrightarrow{[d]}\Gm\longrightarrow 1 .
	\]
	The quotient stack for the weights $\w$ therefore maps naturally to the quotient stack for the weights $\w'$, with relative stabilizer the kernel $\bbm_d$; equivalently, fppf-locally on $\mathcal P_{\w'}$ the category of liftings to $\mathcal P_\w$ is a torsor under $B\bbm_d$.  Hence $\pi$ is a $\bbm_d$-gerbe.

	The induced map on coarse moduli spaces is an isomorphism because the graded rings $k[x_0,\dots,x_n]$ with $\deg(x_i)=w_i$ and with $\deg(x_i)=w_i'$ have the same homogeneous prime ideals after rescaling all degrees by $d$; therefore $\Proj k[x_0,\dots,x_n]_{\w}\cong\Proj k[x_0,\dots,x_n]_{\w'}$.

	Finally, let $x\in U(\overline{k})$ have support $S$.  Its stabilizer consists of $\lambda\in\Gm$ with $\lambda^{w_i}=1$ for all $i\in S$, namely $\bigcap_{i\in S}\bbm_{w_i}=\bbm_{k_S}$.  For the reduced weights $w_i'=w_i/d$ the same computation gives $\bbm_{\gcd\{w_i':i\in S\}}=\bbm_{k_S/d}$.  Hence passing from $\w$ to $\w'$ removes the common subgroup $\bbm_d$ from the stabilizers.
\end{proof}

The reduction $\w\mapsto\w'$ thus exhibits $\cP_\w$ as a gerbe banded by $\bbm_d$ over $\cP_{\w'}$.  Over $\F_q$ the orbit count is, by \cref{thm:stack-classes}, assembled from the local twist groups $H^1(\F_q,\bbm_{k_S})$; we now read the change of these groups under reduction off the $\F_q$-cohomology of the gerbe.  Assume first that $p\nmid d$, so that $\bbm_d$ is \'etale; the wild case is recorded in \cref{rem:wild} below.

\begin{prop}\label{prop:gerbe-cohomology}
Suppose $p\nmid d$.  For each nonempty $S\subseteq T$ the inclusion $\bbm_{k_S/d}\hookrightarrow\bbm_{k_S}$ has quotient $\bbm_d$, giving a short exact sequence of \'etale group schemes over $\F_q$
\[
        1\longrightarrow \bbm_{k_S/d}\longrightarrow \bbm_{k_S}\xrightarrow{(\,\cdot\,)^{k_S/d}}\bbm_d\longrightarrow 1 .
\]
Since $\operatorname{cd}(\F_q)=1$, its long exact $\F_q$-cohomology sequence terminates as the six-term sequence
\[
        0\to\bbm_{k_S/d}(\F_q)\to\bbm_{k_S}(\F_q)\to\bbm_d(\F_q)
        \xrightarrow{\partial} H^1(\F_q,\bbm_{k_S/d})\to H^1(\F_q,\bbm_{k_S})
        \xrightarrow{\pi_\ast} H^1(\F_q,\bbm_d)\to 0 .
\]
In particular $\pi_\ast$ is surjective.  Moreover the gerbe class of $\pi$ in $H^2(\cP_{\w'},\bbm_d)$ restricts at every $\F_q$-point to $H^2(\F_q,\bbm_d)\subseteq\Br(\F_q)=0$, so $\pi$ is neutral over the rational points.
\end{prop}

\begin{proof}
The sequence of group schemes is the standard presentation of $\bbm_d$ as the quotient $\bbm_{k_S}/\bbm_{k_S/d}$, valid since $k_S/d\mid k_S$.  Over a finite field $\operatorname{cd}(\F_q)=1$, so $H^i(\F_q,M)=0$ for $i\ge2$ and every finite \'etale $\F_q$-group $M$; in particular $H^2(\F_q,\bbm_{k_S/d})=0$, which truncates the long exact sequence to the displayed six terms and forces $\pi_\ast$ surjective.  Neutrality is the vanishing $H^2(\F_q,\bbm_d)=0$, equivalently $\Br(\F_q)=0$.
\end{proof}

\begin{cor}\label{cor:gerbe-neutral}
Suppose $p\nmid d$.  The map $\isoclasses(\cP_\w(\F_q))\to\isoclasses(\cP_{\w'}(\F_q))$ is surjective: every $\F_q$-twist of $\cP_{\w'}$ lifts to $\cP_\w$.  Over a coarse point of support $S$ the number of $\F_q$-twists is multiplied, on passing from $\w'$ to $\w$, by
\[
        \frac{\#H^1(\F_q,\bbm_{k_S})}{\#H^1(\F_q,\bbm_{k_S/d})}
        =\frac{\gcd(k_S,q-1)}{\gcd(k_S/d,q-1)}
        =\frac{\#\im(\pi_\ast)}{\#\im(\partial)},
\]
where $\im(\pi_\ast)=H^1(\F_q,\bbm_d)$ has order $\gcd(d,q-1)$ and $\im(\partial)=\operatorname{coker}\bigl(\bbm_{k_S}(\F_q)\to\bbm_d(\F_q)\bigr)$.  The fibre over a given $\w'$-twist is a torsor under $\ker(\pi_\ast)$.
\end{cor}

\begin{proof}
Surjectivity is surjectivity of $\pi_\ast$ in \cref{prop:gerbe-cohomology} together with the lifting of objects guaranteed by neutrality.  The displayed orders are the alternating-product (Euler-characteristic) relation read off the six-term sequence, using $\#\bbm_k(\F_q)=\#H^1(\F_q,\bbm_k)=\gcd(k,q-1)$ from \cref{lem:kummer}.  The fibre statement is exactness at $H^1(\F_q,\bbm_{k_S})$.
\end{proof}

\begin{prop}\label{prop:normalization-monotone}
With $d=\gcd(w_0,\dots,w_n)$ and $\w'=\w/d$, for every $q$ one has $A_{\w'}(q)\le A_\w(q)$, with equality if $\gcd(d,q-1)=1$.
\end{prop}

\begin{proof}
For every nonempty $S$ one has $k_S(\w)=d\,k_S(\w')$, so $\gcd(k_S(\w'),q-1)\mid\gcd(k_S(\w),q-1)$; the support formula of \cref{lem:orbit-formula} gives the inequality term by term.  If $\gcd(d,q-1)=1$ then $\gcd(d\,k_S(\w'),q-1)=\gcd(k_S(\w'),q-1)$ for all $S$, and equality follows.
\end{proof}

\begin{rem}\label{rem:reduction-numerical}
Summing the per-stratum count of \cref{cor:gerbe-neutral} recovers an explicit closed form.  With $d'=\gcd(d,q-1)$, $H=(\F_q^\times)^d$, and $N'(\mu)=\#\{i:\mu^{w'_i}=1\}$ for $\mu\in\F_q^\times$, one has
\[
        A_\w(q)=d'\,A_{\w'}(q)+(d'-1)-\frac{d'}{q-1}\sum_{\mu\in\F_q^\times\setminus H}q^{N'(\mu)}.
\]
Indeed, by \cref{thm:orbit-count}, $(q-1)A_\w(q)+(q-1)=\sum_{\lambda\in\F_q^\times}q^{N(\lambda)}$; since $w_i=d w'_i$ we have $N(\lambda)=N'(\lambda^d)$, and $\lambda\mapsto\lambda^d$ has kernel of order $d'$ and image $H$, so $\sum_\lambda q^{N(\lambda)}=d'\sum_{\mu\in H}q^{N'(\mu)}$; substituting $\sum_{\mu\in H}=\sum_{\mu\in\F_q^\times}-\sum_{\mu\notin H}$ and using the Burnside form for $\w'$ gives the formula.  The factor $d'=\#H^1(\F_q,\bbm_d)$ is the order of the gerbe twist group, and the correction sum measures the failure of $\partial$ to vanish.
\end{rem}

\begin{rem}\label{rem:wild}
If $p\mid d$ the gerbe is no longer tame and $\bbm_d$ is not \'etale, so \cref{prop:gerbe-cohomology} does not apply directly.  The orbit count is nonetheless unaffected: $\gcd(k_S,q-1)$ depends only on the prime-to-$p$ part of $k_S$ (as $p\nmid q-1$), and the support formula of \cref{lem:orbit-formula} holds for all $\w$ unconditionally.
\end{rem}

\begin{thm}\label{thm:partial-rescaling}
Let $\w=(w_0,w_1,\dots,w_n)$ and $\w'=(w_0,\gamma w_1,\dots,\gamma w_n)$ with $\gamma\ge1$ and $\gcd(w_0,\gamma)=1$.  Then
\[
        A_{\w'}(q)-A_\w(q)
        =\sum_{\emptyset\ne S\subseteq\{1,\dots,n\}}(q-1)^{|S|-1}\bigl(\gcd(\gamma k_S,q-1)-\gcd(k_S,q-1)\bigr),
\]
where $k_S=\gcd\{w_i:i\in S\}$.
\end{thm}

\begin{proof}
Apply \cref{thm:orbit-count} to both weight vectors.  If $0\in S$ then $\gcd(w_0,\gamma w_i:i\in S\setminus\{0\})=\gcd(w_0,w_i:i\in S\setminus\{0\})$ because $\gcd(w_0,\gamma)=1$, so those terms cancel.  If $0\notin S$ the corresponding gcd is multiplied by $\gamma$.  The formula follows.
\end{proof}

\begin{cor}\label{cor:partial-equality}
With the notation of \cref{thm:partial-rescaling}, $A_{\w'}(q)=A_\w(q)$ if $\gcd(\gamma k_S,q-1)=\gcd(k_S,q-1)$ for every nonempty $S\subseteq\{1,\dots,n\}$; in particular if $\gcd(\gamma,q-1)=1$.
\end{cor}

\subsection{Numerical examples}\label{subsec:normalization-numerical}

The coarse spaces attached to $\w$ and $\w/d$ are isomorphic, but the orbit count may change because the Kummer groups $H^1(\F_q,\bbm_{k_S})$ change.

\begin{exa}[A nontrivial reduction effect]\label{exa:normalization-24610}
Let $\w=(2,4,6,10)$ and $\w'=(1,2,3,5)$, so $d=2$.  For odd $q$, $\gcd(d,q-1)=2$, and the $\bbm_2$-gerbe contributes nontrivially.  The counts in \cref{tab:normalization-24610} verify \cref{rem:reduction-numerical} and show that the reduction $\w\mapsto\w'$ changes the orbit count although it does not change the coarse variety.

\begin{table}[ht]
\centering
\caption{Orbit counts before and after reduction for $\w=(2,4,6,10)$ and $\w'=(1,2,3,5)$.}
\label{tab:normalization-24610}
\begin{tabular}{c c c c}
\toprule
$q$ & $\gcd(2,q-1)$ & $A_{\w'}(q)$ & $A_{\w}(q)$\\
\midrule
$3$  & $2$ & $41$   & $80$\\
$5$  & $2$ & $157$  & $314$\\
$7$  & $2$ & $403$  & $804$\\
$11$ & $2$ & $1469$ & $2936$\\
\bottomrule
\end{tabular}
\end{table}
\end{exa}

\begin{exa}\label{exa:partial-rescaling-table}
Let $\w=(1,2,3,5)$ and $\w'=(1,4,6,10)$, the case $w_0=1$, $\gamma=2$ of \cref{thm:partial-rescaling}.  The differences in \cref{tab:partial-rescaling} agree with that formula.

\begin{table}[ht]
\centering
\caption{Comparison of orbit counts for $\w=(1,2,3,5)$ and $\w'=(1,4,6,10)$.}
\label{tab:partial-rescaling}
\begin{tabular}{c c c c}
\toprule
$q$ & $A_{\w'}(q)$ & $A_{\w}(q)$ & $A_{\w'}(q)-A_{\w}(q)$\\
\midrule
$5$  & $189$  & $157$  & $32$\\
$7$  & $461$  & $403$  & $58$\\
$11$ & $1605$ & $1469$ & $136$\\
\bottomrule
\end{tabular}
\end{table}
\end{exa}

\section{Singular and smooth strata}\label{sec:strata}

Let $d=\gcd(w_0,\dots,w_n)$.  The generic stabilizer of $\cP_\w$ is $\bbm_d$, and a support stratum has stabilizer $\bbm_{k_S}$.  The non-generic isotropy strata are those with $k_S>d$; on the coarse variety these are the singular strata in the usual well-formed sense \cite{Mori1975,Hosgood}.

\begin{defn}
Define
\[
        A_\w^{\sm}(q)=\sum_{\substack{\emptyset\ne S\subseteq T\\ k_S=d}}(q-1)^{|S|-1}\gcd(k_S,q-1),
        \qquad
        A_\w^{\sing}(q)=\sum_{\substack{\emptyset\ne S\subseteq T\\ k_S>d}}(q-1)^{|S|-1}\gcd(k_S,q-1).
\]
\end{defn}

\begin{thm}\label{thm:strata-decomposition}
$A_\w(q)=A_\w^{\sm}(q)+A_\w^{\sing}(q)$.
\end{thm}

\begin{proof}
The conditions $k_S=d$ and $k_S>d$ partition the nonempty supports, so this is the support formula of \cref{thm:orbit-count} split accordingly.
\end{proof}

\begin{cor}\label{cor:trivial-twist}
If $\gcd(w_i,q-1)=1$ for all $i$, then $\gcd(k_S,q-1)=1$ for all nonempty $S$, and $A_\w(q)=\#\bP^n_\w(\F_q)=(q^{n+1}-1)/(q-1)$.
\end{cor}

\begin{proof}
Every $k_S$ is then coprime to $q-1$, so the formula of \cref{thm:orbit-count} reduces to the coarse count of \cref{cor:coarse-count}.
\end{proof}

\begin{exa}\label{exa:weights-1235}
Let $\w=(1,2,3,5)$, so $d=1$.  The non-generic one-dimensional support strata have stabilizers $\bbm_2,\bbm_3,\bbm_5$, so
\[
        A_\w^{\sing}(q)=\gcd(2,q-1)+\gcd(3,q-1)+\gcd(5,q-1),
        \qquad
        A_\w^{\sm}(q)=q^3+q^2+q-2,
\]
and $A_\w(q)=q^3+q^2+q-2+\gcd(2,q-1)+\gcd(3,q-1)+\gcd(5,q-1)$.
\end{exa}

\subsection{Computational examples}\label{subsec:strata-computational}

The following formulas show how the support formula separates the smooth- and singular-isotropy parts; the corresponding values for $q=p^a$ with $p\in\{3,5,7,11\}$ and $1\le a\le5$ are collected in Appendix~\ref{app:tables}.

\begin{exa}\label{exa:strata-24610}
For $\w=(2,4,6,10)$, $d=2$:
\[
\begin{split}
        A_\w^{\sing}(q)	& 	=\gcd(4,q-1)+\gcd(6,q-1)+\gcd(10,q-1),     \\
        A_\w^{\sm}(q)	&	=(q^3+q^2+q-2)\gcd(2,q-1).
\end{split}
\]
The values are recorded in \cref{tab:appendix-24610}.
\end{exa}

\begin{exa}\label{exa:strata-1235}
For $\w=(1,2,3,5)$, $A_\w^{\sing}(q)=\gcd(2,q-1)+\gcd(3,q-1)+\gcd(5,q-1)$ and $A_\w^{\sm}(q)=q^3+q^2+q-2$.  The values are recorded in \cref{tab:appendix-1235}.
\end{exa}

\begin{exa}\label{exa:strata-161421}
For $\w=(1,6,14,21)$, the singular-isotropy supports occur in dimensions zero and one:
\[
\begin{aligned}
        A_\w^{\sing}(q)
        ={}& \gcd(6,q-1)+\gcd(14,q-1)+\gcd(21,q-1)\\
        & +(q-1)\bigl(\gcd(2,q-1)+\gcd(3,q-1)+\gcd(7,q-1)\bigr),
\end{aligned}
\]
while $A_\w^{\sm}(q)=q^3+q^2-2q+1$.  The values are recorded in \cref{tab:appendix-161421}.
\end{exa}

The singular contribution $A_\w^{\sing}(q^r)$ need not be monotone in $r$: the twist orders $\gcd(k_S,q^r-1)$ fluctuate with the extension degree, so the smooth and singular parts oscillate even as the total grows.

\section{Twist zeta functions}\label{sec:zeta}

We assemble the extension-field orbit counts into a generating function.

\begin{defn}\label{def:zeta}
The twist zeta function of $\cP_\w$ over $\F_q$ is
\[
        \Ztw(\cP_\w,t)=\exp\left(\sum_{r\ge1}\frac{A_\w(q^r)}r t^r\right).
\]
\end{defn}

This records the unweighted isomorphism-class counts of stacky rational forms over finite extensions.  It is not the Hasse--Weil zeta function of the coarse variety,
\(Z(\bP_\w^n,t)=\prod_{j=0}^{n}(1-q^j t)^{-1},\)
nor the stacky mass zeta function, which by \cref{prop:mass-equals-coarse} coincides with $Z(\bP_\w^n,t)$.

For a positive integer $m$ prime to $p$ let $o_m=\ord_m(q)$, set $\ord_1(q)=1$, and let $\mathcal D(k)=\{m:m\mid k,\ p\nmid m\}$.

\begin{thm}\label{thm:zeta-rational}
The function $\Ztw(\cP_\w,t)$ is rational.  More precisely,
\[
\begin{split}
        \Ztw(\cP_\w,t)  	& =        \prod_{\emptyset\ne S\subseteq T}\prod_{m\in\mathcal D(k_S)}\prod_{j=0}^{|S|-1}        \left(1-(q^j t)^{o_m}\right)^{-c(S,m,j)},  \\
        c(S,m,j)	& =\frac{\varphi(m)}{o_m}\binom{|S|-1}{j}(-1)^{|S|-1-j},
\end{split}        
\]
and the exponents $c(S,m,j)$ are integers.
\end{thm}

\begin{proof}
By \cref{cor:extensions}, $A_\w(q^r)=\sum_{\emptyset\ne S\subseteq T}(q^r-1)^{|S|-1}\gcd(k_S,q^r-1)$.  For $k\ge1$,
\[
        \gcd(k,q^r-1)=\sum_{\substack{m\mid k\\ m\mid q^r-1}}\varphi(m)=\sum_{m\in\mathcal D(k)}\varphi(m)\,\I(o_m\mid r),
\]
and $(q^r-1)^{|S|-1}=\sum_{j=0}^{|S|-1}\binom{|S|-1}{j}(-1)^{|S|-1-j}q^{rj}$.  Thus $A_\w(q^r)$ is a finite sum of terms $C(S,m,j)q^{rj}\I(o_m\mid r)$ with $C(S,m,j)=\varphi(m)\binom{|S|-1}{j}(-1)^{|S|-1-j}$.  Since
\[
        \sum_{r\ge1}\frac{q^{rj}t^r}{r}\I(o_m\mid r)
        =\sum_{a\ge1}\frac{(q^jt)^{o_m a}}{o_m a}
        =-\frac1{o_m}\log\bigl(1-(q^jt)^{o_m}\bigr),
\]
exponentiating yields the displayed product with exponent $-C(S,m,j)/o_m=-c(S,m,j)$.  As $o_m\mid\varphi(m)$ for $m$ prime to $p$, the $c(S,m,j)$ are integers.
\end{proof}

\begin{prop}\label{prop:zeta-strata}
With $\Ztw^{\sm}(\cP_\w,t)=\exp(\sum_{r\ge1}A_\w^{\sm}(q^r)t^r/r)$ and $\Ztw^{\sing}(\cP_\w,t)=\exp(\sum_{r\ge1}A_\w^{\sing}(q^r)t^r/r)$, both factors are rational and
\[
        \Ztw(\cP_\w,t)=\Ztw^{\sm}(\cP_\w,t)\,\Ztw^{\sing}(\cP_\w,t).
\]
\end{prop}

\begin{proof}
Rationality follows from the proof of \cref{thm:zeta-rational} after restricting the product to supports with $k_S=d$ or $k_S>d$; the factorization follows from $A_\w(q^r)=A_\w^{\sm}(q^r)+A_\w^{\sing}(q^r)$ for all $r$.
\end{proof}

\begin{prop}\label{prop:zeros-poles}
Every reciprocal zero or pole of $\Ztw(\cP_\w,t)$ has the form $\alpha=\zeta q^j$ with $0\le j\le n$ and $\zeta$ a root of unity; in particular it is an algebraic integer of absolute value $q^j$.
\end{prop}

\begin{proof}
Immediate from the product of \cref{thm:zeta-rational}, whose factors are $1-(q^j t)^{o_m}$.
\end{proof}


\begin{exa}\label{exa:zeta-1235}
For $\w=(1,2,3,5)$, by \cref{exa:weights-1235},
\[
        A_\w(q^r)=q^{3r}+q^{2r}+q^r-2+\gcd(2,q^r-1)+\gcd(3,q^r-1)+\gcd(5,q^r-1).
\]
The zeta function is the coarse-type contribution from $q^{3r}+q^{2r}+q^r-2$ times three Kummer factors: for a prime $\ell\in\{2,3,5\}$ with $o_\ell=\ord_\ell(q)$,
\[
        \exp\left(\sum_{r\ge1}\frac{\gcd(\ell,q^r-1)}r t^r\right)=(1-t)^{-1}\left(1-t^{o_\ell}\right)^{-(\ell-1)/o_\ell}.
\]
\end{exa}

\begin{exa}[Comparison with the reduced weights]\label{exa:zeta-24610}
Let $\w=(2,4,6,10)$ and $\w'=(1,2,3,5)$.  Although $\bP^3_\w\cong\bP^3_{\w'}$ as coarse varieties, the twist zeta functions differ.  Over $\F_{q^r}$,
\[
        A_\w^{\sm}(q^r)=\gcd(2,q^r-1)\bigl(q^{3r}+q^{2r}+q^r-2\bigr),
        \qquad
        A_{\w'}^{\sm}(q^r)=q^{3r}+q^{2r}+q^r-2,
\]
and
\[
        A_\w^{\sing}(q^r)=\gcd(4,q^r-1)+\gcd(6,q^r-1)+\gcd(10,q^r-1),
\]
whereas $A_{\w'}^{\sing}(q^r)=\gcd(2,q^r-1)+\gcd(3,q^r-1)+\gcd(5,q^r-1)$.  Hence $\Ztw(\cP_\w,t)$ and $\Ztw(\cP_{\w'},t)$ are structurally different---the zeta-level reflection of the $\bbm_2$-gerbe of \cref{subsec:normalization-numerical}.
\end{exa}

\section{Functional equations and the multiplicity spectrum}\label{sec:fe}

The rationality of $\Ztw(\cP_\w,t)$ established in \cref{sec:zeta} raises the question of its symmetries: does the twist zeta function satisfy a functional equation, as the Hasse--Weil zeta function of a smooth proper variety does?  We answer this in two stages.  First, the product of \cref{thm:zeta-rational} is naturally indexed by support: writing
\[
\begin{split}
        Z_S(t)		&	=\prod_{m\in\mathcal D(k_S)}\prod_{j=0}^{|S|-1}\bigl(1-(q^jt)^{o_m}\bigr)^{-c(S,m,j)}, \\
        \Ztw(\cP_\w,t)	&	=\prod_{\emptyset\ne S\subseteq T}Z_S(t),
\end{split}
\]
each factor $Z_S$ obeys a functional equation centred at the dimension $\delta_S=|S|-1$ of its stratum, driven by a symmetry of the exponents under $j\mapsto\delta_S-j$.  Second, a sector decomposition of the orbit count determines exactly when these per-stratum symmetries align into a single global functional equation.

\begin{lem}\label{lem:exponent-symmetry}
Let $\delta=|S|-1$.  For all $m\in\mathcal D(k_S)$ and $0\le j\le\delta$,
\[
        c(S,m,\delta-j)=(-1)^{\delta}\,c(S,m,j).
\]
\end{lem}

\begin{proof}
From
\[
c(S,m,j)=\tfrac{\varphi(m)}{o_m}\binom{\delta}{j}(-1)^{\delta-j}
\]
and $\binom{\delta}{\delta-j}=\binom{\delta}{j}$ we get $c(S,m,\delta-j)=\tfrac{\varphi(m)}{o_m}\binom{\delta}{j}(-1)^{j}$, and $(-1)^{j}=(-1)^{\delta}(-1)^{\delta-j}$.
\end{proof}

Write $k_S^{(p')}$ for the prime-to-$p$ part of $k_S$, and note $\sum_{m\in\mathcal D(k_S)}\varphi(m)=k_S^{(p')}$.

\begin{thm}\label{thm:fe}
Let $\delta=|S|-1$.  The multiset of reciprocal zeros and poles of $Z_S(t)$, counted with multiplicity, is invariant under $\alpha\mapsto q^{\delta}/\alpha$; this involution preserves zeros and poles when $\delta$ is even and interchanges them when $\delta$ is odd.  Explicitly,
\[
        Z_S\!\left(\frac{1}{q^{\delta}t}\right)=\gamma_S\,t^{\sigma_S}\,Z_S(t)^{(-1)^{\delta}},
\]
where
\[
        \sigma_S=\begin{cases} k_S^{(p')},&\delta=0,\\ 0,&\delta\ge1,\end{cases}
        \qquad
        \gamma_S=\begin{cases}(-1)^{\nu_S},&\delta=0,\ \ \nu_S=\displaystyle\sum_{m\in\mathcal D(k_S)}\tfrac{\varphi(m)}{o_m},\\[1mm]
        q^{-k_S^{(p')}},&\delta=1,\\ 1,&\delta\ge2.\end{cases}
\]
\end{thm}

\begin{proof}
The substitution $t\mapsto 1/(q^{\delta}t)$ gives
\[
        1-(q^j t)^{o_m}\ \longmapsto\ 1-q^{-(\delta-j)o_m}t^{-o_m}
        =-\,q^{-(\delta-j)o_m}t^{-o_m}\bigl(1-(q^{\delta-j}t)^{o_m}\bigr).
\]
Raising to the power $-c(S,m,j)$ and multiplying over all $(m,j)$,
\[
        Z_S\!\left(\frac{1}{q^{\delta}t}\right)
        =\Bigl[\prod_{m,j}(-1)^{-c(S,m,j)}q^{(\delta-j)o_m c(S,m,j)}t^{o_m c(S,m,j)}\Bigr]
        \cdot\prod_{m,j}\bigl(1-(q^{\delta-j}t)^{o_m}\bigr)^{-c(S,m,j)}.
\]
In the last product reindex $j\mapsto\delta-j$ and apply \cref{lem:exponent-symmetry}: the exponent becomes $-(-1)^{\delta}c(S,m,j)$, so the product equals $Z_S(t)^{(-1)^{\delta}}$.  This proves the displayed shape, and the root-multiset statement is the same identity read off the factors $1-(q^{j}t)^{o_m}$.

For the constant, combine $\sum_{m}\varphi(m)=k_S^{(p')}$ with the binomial identities
\[
        \sum_{j=0}^{\delta}\binom{\delta}{j}(-1)^{\delta-j}=\I(\delta=0),
        \qquad
        \sum_{j=0}^{\delta}j\binom{\delta}{j}(-1)^{\delta-j}=\I(\delta=1),
\]
the values at $x=1$ of $(x-1)^{\delta}$ and of $x\,\tfrac{d}{dx}(x-1)^{\delta}$.  The $t$-exponent is $\sum_{m,j}o_m c(S,m,j)=k_S^{(p')}\,\I(\delta=0)$, giving $\sigma_S$; the $q$-exponent $\sum_{m,j}(\delta-j)o_m c(S,m,j)$ equals $0$, $-k_S^{(p')}$, $0$ for $\delta=0$, $\delta=1$, $\delta\ge2$ respectively; and the sign has parity $\sum_{m,j}c(S,m,j)=\nu_S\,\I(\delta=0)$.
\end{proof}

Grouping by dimension, $\Ztw(\cP_\w,t)=\prod_{\delta=0}^{n}Z^{(\delta)}(t)$ with $Z^{(\delta)}=\prod_{|S|-1=\delta}Z_S$, and each $Z^{(\delta)}$ inherits from \cref{thm:fe} a functional equation centred at $q^{\delta}$.  A single functional equation for $\Ztw$ requires all contributing factors to share one centre; the obstruction is measured by the multiplicity spectrum.

\begin{defn}\label{def:spectrum}
Fix $q$ and write, by \cref{thm:zeta-rational},
\[
\Ztw(\cP_\w,t)=\prod_{\alpha}(1-\alpha t)^{-b(\alpha)}
\]
over reciprocal roots $\alpha=\zeta q^{j}$.  The \emph{multiplicity spectrum} is $\mathbf b_\w(q)=(b_0,\dots,b_n)$, where $b_j=\sum_{|\alpha|=q^{j}}b(\alpha)$ is the total signed multiplicity at radius $q^{j}$.  In the split regime $q\equiv1\pmod{\operatorname{lcm}_i w_i}$ one has $o_m=1$ throughout and
\[
\begin{split}
b_j		&	=\sum_{\substack{\emptyset\ne S\subseteq T\\ |S|-1\ge j}}\binom{|S|-1}{j}(-1)^{|S|-1-j}\gcd(k_S,q-1),\\
 A_\w(q^r)	&	=\sum_{j=0}^{n}b_j\,q^{rj}.
\end{split}
\]
\end{defn}

\begin{thm}\label{thm:fe-global}
The twist zeta function $\Ztw(\cP_\w,t)$ satisfies a single functional equation under $t\mapsto 1/(q^{n}t)$, relating it to $\Ztw(\cP_\w,t)^{\pm1}$ up to an explicit monomial, if and only if the multiplicity spectrum is \emph{palindromic}:
\[
        b_j=b_{n-j}\qquad\text{for all }0\le j\le n.
\]
\end{thm}

\begin{proof}
The substitution $t\mapsto1/(q^nt)$ sends a reciprocal root of radius $q^{j}$ to one of radius $q^{n-j}$, so it preserves the root multiset---which is what a functional equation of this centre asserts---iff $b_j=b_{n-j}$ for every $j$.  When this holds, pairing each factor with its image as in \cref{thm:fe} produces the stated identity.
\end{proof}

\subsection{The spectrum via isotropy sectors}\label{subsec:sectors}

By \cref{thm:fe-global}, a single functional equation for $\Ztw(\cP_\w,t)$ is equivalent to palindromy of the multiplicity spectrum.  To decide when this holds we expand the orbit count a second way, grouping scalars in the Burnside formula by their exact order rather than vectors by their support.  For an integer $e\ge1$ set
\[
        n(e)=\#\{\,i\in T : e\mid w_i\,\},
\]
so that $n(1)=n+1$, and $n(e)\ge1$ only for the finitely many $e$ dividing some weight.  The locus where $\bbm_e$ acts trivially is the coordinate subspace on $\{i: e\mid w_i\}$, with coarse space $\bP^{\,n(e)-1}$; we call it the \emph{$e$-sector}.

\begin{thm}\label{thm:sector}
For every $r\ge1$,
\[
A_\w(q^r)        =\sum_{\substack{e\ge1,\ p\nmid e\\ o_e\mid r}}\varphi(e)\,        \frac{q^{r\,n(e)}-1}{q^r-1}        =\sum_{\substack{e\ge1,\ p\nmid e\\ o_e\mid r}}\varphi(e)\,       \#\bP^{\,n(e)-1}(\F_{q^r}),
\]
and consequently
\[
  \Ztw(\cP_\w,t)        =\prod_{\substack{e\ge1,\ p\nmid e\\ n(e)\ge1}}\    \prod_{j=0}^{n(e)-1}\bigl(1-(q^jt)^{o_e}\bigr)^{-\varphi(e)/o_e},
\]
a finite product with positive integer exponents.  In particular $\Ztw(\cP_\w,t)^{-1}\in\ZZ[t]$: the twist zeta function has poles only.
\end{thm}

\begin{proof}
In the Burnside formula of \cref{lem:burnside} over $\F_{q^r}$, group the scalars by exact order: for each $e\mid q^r-1$ (necessarily prime to $p$) there are $\varphi(e)$ scalars of order $e$, and $N(\lambda)=n(e)$ for each, since $\lambda^{w_i}=1$ iff $e\mid w_i$.  Hence
\[
        A_\w(q^r)=\frac1{q^r-1}\sum_{e\mid q^r-1}\varphi(e)\bigl(q^{r\,n(e)}-1\bigr),
\]
which is the first display, since $e\mid q^r-1$ iff $o_e\mid r$ and the terms with $n(e)=0$ vanish.

Expanding $(q^{r\,n(e)}-1)/(q^r-1)=\sum_{j=0}^{n(e)-1}q^{rj}$ and exponentiating as in the proof of \cref{thm:zeta-rational} gives the product, with exponent $\varphi(e)/o_e\in\ZZ$ since $o_e\mid\varphi(e)$.
\end{proof}

\begin{cor}\label{cor:spectrum-sector}
For every $q$ of characteristic $p$,
\[
b_j=\sum_{\substack{e\ge1,\ p\nmid e\\ n(e)\ge j+1}}\varphi(e),        \qquad 0\le j\le n .
\]
In particular:

\begin{enumerate}[label=\textup{(\roman*)}]

\item the spectrum depends only on $\w$ and $p$, not on $q$;

\item $b_0\ge b_1\ge\cdots\ge b_n=d^{(p')}$;

\item $b_j\ge1$ for all $j$, and for reduced $\w$ the defect coefficients are $b_j-1=\sum_{e\ge2,\,p\nmid e,\,n(e)\ge j+1}\varphi(e)\ge0$.
\end{enumerate}
\end{cor}

\begin{proof}
In the product of \cref{thm:sector} the factor attached to $(e,j)$ contributes $o_e$ reciprocal poles of modulus $q^{j}$, each of multiplicity $\varphi(e)/o_e$, and there are no zeros; summing over the sectors with $n(e)\ge j+1$ gives the formula.  The condition $n(e)\ge j+1$ is nested in $j$, whence (ii); the sector $e=1$ contributes $1$ to every $b_j$, whence (iii); and $n(e)=n+1$ iff $e\mid d$, whence $b_n=d^{(p')}$.
\end{proof}

\begin{exa}\label{exa:spectra}
For $\w=(1,2,3,5)$ the sectors with $n(e)\ge1$ are $e=1$, with $n(1)=4$, and $e=2,3,5$, with $n(e)=1$; so $b_0=1+\varphi(2)+\varphi(3)+\varphi(5)=8$ and $\mathbf b_\w=(8,1,1,1)$.  For $\w=(1,1,2,2)$ the only nontrivial sector is $e=2$ with $n(2)=2$, so $\mathbf b_\w=(2,2,1,1)$.  For $\w=(1,6,14,21)$ the sectors $e=2,3,7$ have $n(e)=2$ and $e=6,14,21$ have $n(e)=1$, so $b_1=10$, $b_0=30$, and $\mathbf b_\w=(30,10,1,1)$.  In the split regime these give $A_\w(q)=q^3+q^2+q+8$, $q^3+q^2+2q+2$, and $q^3+q^2+10q+30$ respectively; direct evaluation of the support count of \cref{thm:orbit-count} at $q=31$, $5$, and $43$ confirms the three spectra.  Each is nonincreasing with $b_3=1$, and none is palindromic.
\end{exa}

\begin{thm}\label{thm:palindromy-classification}
The following are equivalent:
\begin{enumerate}[label=\textup{(\alph*)}]
\item $\Ztw(\cP_\w,t)$ satisfies a single functional equation under $t\mapsto1/(q^nt)$;
\item the multiplicity spectrum is palindromic;
\item $b_0=b_n$;
\item the prime-to-$p$ parts $w_0^{(p')},\dots,w_n^{(p')}$ are all equal.
\end{enumerate}
When these hold, with common value $c=w_i^{(p')}$, the spectrum is $(c,\dots,c)$.  For reduced $\w$ they hold if and only if every $w_i$ is a power of $p$, in which case all twists are trivial (\cref{cor:trivial-twist}) and $\Ztw(\cP_\w,t)=Z(\bP^n,t)$.
\end{thm}

\begin{proof}
(a)$\Leftrightarrow$(b) is \cref{thm:fe-global}, and (b)$\Rightarrow$(c) is immediate.  Since the spectrum is nonincreasing by \cref{cor:spectrum-sector}(ii), $b_0=b_n$ forces it constant, giving (c)$\Rightarrow$(b).

For (c)$\Leftrightarrow$(d): by \cref{cor:spectrum-sector}, $b_0$ sums $\varphi(e)$ over the prime-to-$p$ $e$ dividing some weight, and $b_n$ over those dividing all weights.  As $\varphi(e)>0$ and the second index set is contained in the first, $b_0=b_n$ iff the two sets coincide; taking $e=w_i^{(p')}$ then gives $w_i^{(p')}\mid w_j^{(p')}$ for all $i,j$, so all prime-to-$p$ parts are equal, and the converse is clear.  If (d) holds with common value $c$, then $n(e)=n+1$ for $e\mid c$ and $n(e)=0$ for every other prime-to-$p$ $e$, so $b_j=\sum_{e\mid c}\varphi(e)=c$ for all $j$.

For reduced $\w$, writing $w_i=c\,p^{a_i}$ gives $1=d=c\,p^{\min a_i}$, so $c=1$ and each $w_i$ is a $p$-power; then every $k_S$ is a $p$-power, $\gcd(k_S,q^r-1)=1$ for all $r$, and $A_\w(q^r)=\#\bP^n(\F_{q^r})$.
\end{proof}

\begin{cor}\label{cor:no-self-dual}
Let $\w$ be reduced.  If some integer $e\ge2$ prime to $p$ divides one of the weights---that is, if $\cP_\w$ has a nontrivial tame isotropy sector---then $\Ztw(\cP_\w,t)$ admits no single global functional equation, for any $q$.  The orders $o_e$ redistribute the reciprocal poles $\zeta q^{j}$ around the circles $|\alpha|=q^{j}$ but never between them, so the spectrum, and with it the failure of palindromy, is the same in every regime.
\end{cor}

\begin{exa}\label{exa:fe-1235}
For $\w=(1,2,3,5)$, \cref{thm:sector} gives
\[
        \Ztw(\cP_\w,t)        =\frac{1}{(1-t)(1-qt)(1-q^2t)(1-q^3t)}        \prod_{\ell\in\{2,3,5\}}\bigl(1-t^{o_\ell}\bigr)^{-(\ell-1)/o_\ell}.
\]
The first factor is the $e=1$ sector, the Hasse--Weil zeta function $Z(\bP^3,t)$ of the coarse space, self-dual under $t\mapsto1/(q^3t)$ up to a monomial; each prime $\ell$ contributes a dimension-$0$ sector centred at $q^0=1$.  The centres do not align, and indeed the spectrum $(8,1,1,1)$ of \cref{exa:spectra} is not palindromic: the prime-to-$p$ parts of the weights are not all equal.

The example also displays the regime-independence of the spectrum.  In the split regime each sector factor is $(1-t)^{-(\ell-1)}$, a pole at $\alpha=1$ of multiplicity $\ell-1$; for $q\equiv-1\pmod3$ the $3$-sector becomes $(1-t^{2})^{-1}$, and the multiplicity $2$ splits between $\alpha=1$ and $\alpha=-1$.  The poles move around the circle $|\alpha|=1$ but never leave it, and the spectrum is unchanged.
\end{exa}

\begin{rem}\label{rem:perspective}
The classification just proved is a negative result, and its meaning deserves a remark.  For ordinary projective space $\Ztw$ is the Hasse--Weil zeta function and inherits its classical functional equation; one might have hoped that some weighted stacks---perhaps for special $q$, exploiting the arithmetic of the orders $o_e$---recover this symmetry while retaining nontrivial twists.  \cref{thm:palindromy-classification} closes the door: the spectrum is blind to $q$, and the symmetry survives precisely when the twist theory is trivial.  The obstruction and the invariant are one phenomenon.  The isotropy that makes $A_\w$ finer than the coarse count is exactly what tilts the spectrum, so a weighted stack pays for its twists with its functional equation.

The question reopens one level down, where it is no longer combinatorial.  For a subvariety $\mathcal X\subset\cP_\w$ the stratum counts carry Frobenius eigenvalues, and the spectrum of $\Ztw(\mathcal X,t)$ genuinely depends on $q$; the diagonal case is treated in \cref{sec:hypersurfaces}, where the question of a functional equation is posed but not resolved.  Heights, entropy, and the associated zeta functions over finite fields are developed in the companion paper \cite{ss-heights}.
\end{rem}

\section{Weighted complete intersections}\label{sec:hypersurfaces}

The same method applies to weighted subvarieties, stated for the affine cone to avoid scheme-theoretic ambiguity.  Let $F=(f_1,\dots,f_m)$ be weighted homogeneous polynomials in $\F_q[x_0,\dots,x_n]$ for the weights $\w$, let $X=V(F)\subseteq\bP^n_\w$, and let
\(\widehat X=\{x\in\bbA^{n+1}:f_1(x)=\cdots=f_m(x)=0\}\)
be the affine cone.  For a nonempty $S\subseteq T$ put $N_X(S;q)=\#\{x\in\widehat X(\F_q):\Supp(x)=S\}$.

\begin{defn}
	The orbit count of $X$ is $A_X(q)=\#\big((\widehat X(\F_q)\setminus\{0\})/\F_q^\times\big)$.
\end{defn}

\begin{prop}\label{prop:hypersurface-count}
	With the notation above,
	\[
	A_X(q)=\sum_{\emptyset\ne S\subseteq T}\frac{N_X(S;q)}{q-1}\gcd(k_S,q-1).
	\]
\end{prop}

\begin{proof}
	The equations are invariant under the weighted $\F_q^\times$-action.  On cone points of fixed support $S$ every point has stabilizer $\bbm_{k_S}(\F_q)$ of order $\gcd(k_S,q-1)$, so every orbit there has size $(q-1)/\gcd(k_S,q-1)$.  Dividing $N_X(S;q)$ by this size and summing over $S$ gives the formula.
\end{proof}

\begin{rem}
	This is an orbit formula on the rational cone, not the number of rational points of the coarse subvariety; the two agree only when all local Kummer twist factors are trivial over $\F_q$.
\end{rem}

\begin{exa}\label{exa:hypersurface-P12}
	Let $X\subset\bP_{(1,2)}$ be defined by $y=x^2$.  The nonzero rational cone points $(a,a^2)$, $a\in\F_q^\times$, all have support $\{0,1\}$ and form one orbit, so $A_X(q)=1$.
\end{exa}

\begin{rem}\label{rem:qwag}
	In the evaluation codes of weighted algebraic-geometry constructions \cite{shaska2025}, the relevant count of evaluation points is governed by the coarse rational-point count, not by $A_\w(q)$.  A weighted-homogeneous polynomial $f$ of degree $d$ satisfies $f(\mu\cdot v)=\mu^{d}f(v)$ for $\mu\in\F_q^\times$, so its vanishing and the value $f(v)$ up to the scalar $\mu^d$ are constant along an $\F_q^\times$-orbit; the distinct $\F_q^\times$-orbits over a fixed coarse point counted by $\gcd(k_S,q-1)$ are scaled copies of one another and supply no additional evaluation data.  Thus $A_\w(q)$ does not enter the code parameters, consistently with the observation \cite{shaska2025} that weighted ambient spaces give no automatic gain in rational points; $A_\w(q)$ is the right invariant for twisted objects and rational representatives rather than for the coarse evaluation set.
\end{rem}

\subsection{Twist zeta functions of diagonal hypersurfaces}\label{subsec:diagonal}

Let $\mathcal X=[(\widehat X\setminus\{0\})/\Gm]\subseteq\cP_\w$ be the closed substack defined by $F$.  Exactly as in \cref{thm:stack-classes}, isomorphism classes of objects of $\mathcal X(\F_q)$ are $\F_q^\times$-orbits on $\widehat X(\F_q)\setminus\{0\}$, so $A_X(q)=\#\isoclasses(\mathcal X(\F_q))$.

\begin{defn}\label{def:hyp-zeta}
	The twist zeta function of $\mathcal X$ over $\F_q$ is
	\[
	\Ztw(\mathcal{X}, t) = \exp \left( \sum_{r \ge 1} \frac{A_X(q^r)}{r} t^r \right).
	\]
\end{defn}

We compute it for diagonal hypersurfaces.  Consider
\[
F(x_0, \dots , x_n) = a_0 x_0^{m_0} + a_1 x_1^{m_1} + \cdots  + a_n x_n^{m_n},
\]
with $a_i \in \F_q^\times$ and $m_i w_i = d$ for a fixed degree $d$.  Fix a nontrivial additive character $\psi$ of $\F_Q$ and for a multiplicative character $\chi$ of $\F_Q^\times$ set $G_Q(\chi)=\sum_{u\in\F_Q^\times}\chi(u)\psi(u)$, so that $G_Q(\varepsilon)=-1$ for the trivial character $\varepsilon$ and $|G_Q(\chi)|=Q^{1/2}$ for $\chi\ne\varepsilon$.

\begin{lem}\label{lem:diagonal-count}
Let $Q$ be any power of $q$ and $\emptyset\ne S\subseteq T$.  Then
\[
	N_X(S;Q)=\frac{(Q-1)^{|S|}}{Q}+\frac{Q-1}{Q}\sum_{\substack{(\chi_i)_{i\in S}\\ \chi_i^{m_i}=\varepsilon,\ \prod_{i\in S}\chi_i=\varepsilon}}\ \prod_{i\in S}\overline\chi_i(a_i)\,G_Q(\chi_i),
\]
the sum running over all tuples of multiplicative characters of $\F_Q^\times$, trivial components allowed, with $\chi_i^{m_i}=\varepsilon$ and trivial product.
\end{lem}

\begin{proof}
Orthogonality of additive characters gives
\[
	N_X(S;Q)=\frac1Q\sum_{t\in\F_Q}\ \sum_{x\in(\F_Q^\times)^S}\psi\Bigl(t\sum_{i\in S} a_ix_i^{m_i}\Bigr).
\]
The term $t=0$ contributes $(Q-1)^{|S|}/Q$.  For $t\ne0$ and each $i$, substituting $u=x_i^{m_i}$ and using $\#\{x:x^{m_i}=u\}=\sum_{\chi^{m_i}=\varepsilon}\chi(u)$,
\[
	\sum_{x\ne0}\psi(ta_ix^{m_i})=\sum_{\chi^{m_i}=\varepsilon}\ \sum_{u\ne0}\chi(u)\psi(ta_iu)=\sum_{\chi^{m_i}=\varepsilon}\overline\chi(ta_i)\,G_Q(\chi),
\]
where the term $\chi=\varepsilon$ equals $\sum_{u\ne0}\psi(ta_iu)=-1=\overline\varepsilon(ta_i)G_Q(\varepsilon)$.  Multiplying over $i\in S$, expanding, and summing over $t\ne0$, the factor $\sum_{t\ne0}\prod_{i}\overline\chi_i(t)$ equals $Q-1$ when $\prod_i\chi_i=\varepsilon$ and vanishes otherwise.
\end{proof}

Suppose now $q\equiv1\pmod d$, so that every $m_i$ and every $k_S$ divides $q-1$; note this forces $p\nmid d$.  For $\emptyset\ne S\subseteq T$ let $\mathcal C_S$ denote the set of tuples $(\chi_i)_{i\in S}$ of characters of $\F_q^\times$ with $\chi_i^{m_i}=\varepsilon$ and $\prod_{i\in S}\chi_i=\varepsilon$, and for $\chi\in\mathcal C_S$ set
\[
	\omega_{S,\chi}=\prod_{i\in S}\bigl(-\,\overline\chi_i(a_i)\,G_q(\chi_i)\bigr),
	\qquad
	|\omega_{S,\chi}|=q^{\nu(\chi)/2},
\]
where $\nu(\chi)=\#\{i\in S:\chi_i\ne\varepsilon\}$.

\begin{thm}\label{thm:diagonal-rational}
Let $X=V(F)\subset\bP^n_\w$ be the diagonal hypersurface above and $q\equiv1\pmod d$.  Then
\[
	\Ztw(\mathcal X,t)=\prod_{\emptyset\ne S\subseteq T}\ \prod_{j=0}^{|S|-1}\bigl(1-q^{j-1}t\bigr)^{-(-1)^{|S|-1-j}\binom{|S|-1}{j}k_S}\ \prod_{\chi\in\mathcal C_S}\bigl(1-q^{-1}\omega_{S,\chi}\,t\bigr)^{-(-1)^{|S|}k_S}.
\]
In particular $\Ztw(\mathcal X,t)$ is rational, and its reciprocal zeros and poles lie among the numbers $q^{j-1}$ with $0\le j\le n$ and $q^{-1}\omega_{S,\chi}$ of absolute value $q^{\nu(\chi)/2-1}$.
\end{thm}

\begin{proof}
Since $q^r\equiv1\pmod d$ for every $r$, the characters of $\F_{q^r}^\times$ with $\chi^{m_i}=\varepsilon$ are exactly the norm-lifts $\chi\circ\operatorname{N}_{\F_{q^r}/\F_q}$ of the characters of $\F_q^\times$ with the same property: the lift is injective, preserves orders and products, and both sets have $m_i$ elements.  Hence \cref{lem:diagonal-count} over $\F_{q^r}$ is indexed by $\mathcal C_S$ for every $r$.

For $a\in\F_q^\times$ one has $\operatorname{N}_{\F_{q^r}/\F_q}(a)=a^r$, and the Hasse--Davenport relation gives $G_{q^r}(\chi\circ\operatorname{N})=(-1)^{r-1}G_q(\chi)^r$, valid also for $\chi=\varepsilon$ since $(-1)^{r-1}(-1)^r=-1$.  Therefore
\[
	\prod_{i\in S}\overline{\chi_i\circ\operatorname{N}}(a_i)\,G_{q^r}(\chi_i\circ\operatorname{N})
	=(-1)^{(r-1)|S|}\Bigl(\prod_{i\in S}\overline\chi_i(a_i)G_q(\chi_i)\Bigr)^{r}
	=(-1)^{|S|}\,\omega_{S,\chi}^{\,r}.
\]
Combining with \cref{prop:hypersurface-count}, $\gcd(k_S,q^r-1)=k_S$, and the binomial expansion of $(q^r-1)^{|S|-1}$,
\[
	A_X(q^r)=\sum_{\emptyset\ne S\subseteq T}k_S\Biggl[\ \sum_{j=0}^{|S|-1}\binom{|S|-1}{j}(-1)^{|S|-1-j}\,q^{r(j-1)}\ +\ (-1)^{|S|}\sum_{\chi\in\mathcal C_S}\bigl(q^{-1}\omega_{S,\chi}\bigr)^{r}\Biggr].
\]
Every term has the form $c\,\gamma^r$ with $c\in\ZZ$, so summing $A_X(q^r)t^r/r$ over $r$ and exponentiating via $\sum_{r\ge1}\gamma^rt^r/r=-\log(1-\gamma t)$ yields the displayed product.  The absolute values follow from $|G_q(\chi_i)|=q^{1/2}$ for $\chi_i\ne\varepsilon$ and $|G_q(\varepsilon)|=1$.
\end{proof}

\begin{rem}\label{rem:diagonal-checks}
Two internal checks.  For a singleton $S=\{i\}$ the set $\mathcal C_S$ contains only the trivial tuple, with $\omega=1$, and the two factors attached to $S$ cancel---as they must, since $a_ix^{m_i}=0$ has no solution with $x\ne0$ and the stratum contributes nothing.  At the opposite extreme, the full-support tuples with all components nontrivial have $|q^{-1}\omega_{T,\chi}|=q^{(n-1)/2}$, the weight of the middle cohomology of a smooth hypersurface of dimension $n-1$: the twist zeta function sees the primitive Frobenius eigenvalues of $X$, rescaled by the orbit count.
\end{rem}

\begin{exa}\label{exa:P123}
Let $\w=(1,2,3)$ and $X: x_0^6+x_1^3+x_2^2=0$ over $\F_7$, so $d=6$ and $7\equiv1\pmod 6$.  Direct enumeration of the cone gives $N_X(\{0,1\};7)=18$, $N_X(\{1,2\};7)=6$, and $N_X(S;7)=0$ for every other support, whence $A_X(7)=\tfrac{18}{6}+\tfrac{6}{6}=4$ by \cref{prop:hypersurface-count}.  \cref{lem:diagonal-count} reproduces each value: for $S=\{0,1\}$ the admissible tuples are the trivial one and the two with $\chi_0=\overline\chi_1$ cubic, giving $\tfrac{36}{7}+\tfrac{6}{7}\,(1+2\cdot7)=18$; for $S=\{1,2\}$ only the trivial tuple satisfies $\chi_1=\chi_2$ with $\chi_1^3=\chi_2^2=\varepsilon$, giving $\tfrac{36}{7}+\tfrac{6}{7}=6$; for $S=\{0,2\}$ the quadratic character $\rho$ contributes $G_7(\rho)^2=\rho(-1)\,7=-7$, giving $\tfrac{36}{7}+\tfrac{6}{7}(1-7)=0$.
\end{exa}

\begin{rem}\label{rem:diagonal-general}
For general $q$ the count of \cref{lem:diagonal-count} remains valid over every $\F_{q^r}$, but the index set varies with $r$: a tuple whose components have exact orders $e_i$ appears precisely when $o_{\operatorname{lcm}(e_i)}\mid r$, is defined over the corresponding subextension, and the Hasse--Davenport relation applies relative to that base field, attaching the eigenvalues to Galois orbits of tuples.  We do not pursue this bookkeeping here.  The structural point is the contrast with \cref{sec:fe}: the eigenvalues $\omega_{S,\chi}$ depend on $q$ through the Gauss sums, so the multiplicity spectrum of $\Ztw(\mathcal X,t)$ is no longer the $q$-independent combinatorial datum of \cref{cor:spectrum-sector}, and the question of a functional equation for $\Ztw(\mathcal X,t)$ is a genuine arithmetic question, open even for diagonal hypersurfaces.
\end{rem}

\subsection{Intersections of two diagonal hypersurfaces}\label{subsec:two-diagonal}

The method extends to intersections of two diagonal hypersurfaces of the same degree.  Let
\[
F_1=\sum_{i=0}^{n}a_ix_i^{m_i},\qquad F_2=\sum_{i=0}^{n}b_ix_i^{m_i},
\]
with $m_iw_i=d$ for a fixed degree $d$, coefficients $a_i,b_i\in\F_q$ with $(a_i,b_i)\ne(0,0)$, and the points $(a_i:b_i)\in\bP^1(\F_q)$ pairwise distinct; the last condition is the nondegeneracy of the pair.  Let $X=V(F_1,F_2)\subseteq\bP^n_\w$ and $\mathcal X\subseteq\cP_\w$ the associated stack.  Equal degrees are essential: for $t=(t_1:t_2)\in\bP^1$ the pencil member $t_1F_1+t_2F_2$ is again diagonal, with coefficients $\gamma_i(t)=t_1a_i+t_2b_i$, and the count fibres over this pencil.  The zero of $\gamma_i$ is the point $P_i=(b_i:-a_i)$, the unique member of the pencil in which the variable $x_i$ does not appear; by nondegeneracy the $P_i$ are pairwise distinct.

Fix for each $P\in\bP^1(\F_Q)$ a representative in $\F_Q^2\smallsetminus\{0\}$, and for a tuple $\chi=(\chi_i)_{i\in S}$ of multiplicative characters of $\F_Q^\times$ with $\prod_{i\in S}\chi_i=\varepsilon$ set
\[
\begin{split}
J_S(\chi;Q)		&	=\sum_{\substack{P\in\bP^1(\F_Q)\\ \gamma_i(P)\ne0\ \forall i\in S}}\ \prod_{i\in S}\overline\chi_i\bigl(\gamma_i(P)\bigr), \\
\Sigma_{S,i_0}(Q)	&	=\sum_{\substack{(\chi_i)_{i\in S\smallsetminus i_0}\\ \chi_i^{m_i}=\varepsilon,\ \prod\chi_i=\varepsilon}}\ \prod_{i\in S\smallsetminus i_0}\overline\chi_i\bigl(\gamma_i(P_{i_0})\bigr)G_Q(\chi_i);
\end{split}
\]
both are independent of the chosen representatives, since the character products are trivial.

\begin{lem}\label{lem:pencil-count}
Let $Q$ be any power of $q$ and $\emptyset\ne S\subseteq T$.  With $\mathcal C_S$ the set of tuples $(\chi_i)_{i\in S}$ of characters of $\F_Q^\times$ satisfying $\chi_i^{m_i}=\varepsilon$ and $\prod_{i\in S}\chi_i=\varepsilon$,
\[
Q^2\,N_X(S;Q)=(Q-1)^{|S|}+(Q-1)\sum_{\chi\in\mathcal C_S}\Bigl(\prod_{i\in S}G_Q(\chi_i)\Bigr)J_S(\chi;Q)+(Q-1)^2\sum_{i_0\in S}\Sigma_{S,i_0}(Q).
\]
\end{lem}

\begin{proof}
Orthogonality of additive characters applied to both equations gives
\[
Q^2\,N_X(S;Q)=\sum_{t_1,t_2\in\F_Q}\ \prod_{i\in S}T_i(t_1,t_2),\qquad T_i(t_1,t_2)=\sum_{x\in\F_Q^\times}\psi\bigl((t_1a_i+t_2b_i)\,x^{m_i}\bigr).
\]
The term $(t_1,t_2)=(0,0)$ contributes $(Q-1)^{|S|}$.  A nonzero $(t_1,t_2)$ is $\lambda\cdot s(P)$ for unique $\lambda\in\F_Q^\times$ and $P\in\bP^1(\F_Q)$, and as in the proof of \cref{lem:diagonal-count}, $T_i=\sum_{\chi^{m_i}=\varepsilon}\overline\chi(\lambda\gamma_i(P))G_Q(\chi)$ when $\gamma_i(P)\ne0$, while $T_{i}=Q-1$ when $\gamma_{i}(P)=0$.  By nondegeneracy at most one $\gamma_i$ vanishes at $P$.  If none does, expanding the product and summing over $\lambda$ kills all tuples with $\prod\chi_i\ne\varepsilon$ and leaves $(Q-1)\sum_{\chi\in\mathcal C_S}\prod_i\overline\chi_i(\gamma_i(P))G_Q(\chi_i)$; summing over such $P$ gives the middle term.  If $\gamma_{i_0}(P)=0$, i.e.\ $P=P_{i_0}$, the same computation on $S\smallsetminus\{i_0\}$ gives $(Q-1)^2\,\Sigma_{S,i_0}(Q)$.
\end{proof}

Note the recursion: $\Sigma_{S,i_0}$ is precisely the character sum of \cref{lem:diagonal-count} for the single diagonal hypersurface $\sum_{i\ne i_0}\gamma_i(P_{i_0})x_i^{m_i}$, the member of the pencil that forgets the variable $x_{i_0}$.

Suppose now $q\equiv1\pmod d$, so every $m_i$ and every $k_S$ divides $q^r-1$ for all $r$, and $\mathcal C_S$ over $\F_{q^r}$ consists of the norm-lifts of $\mathcal C_S$ over $\F_q$, as in the proof of \cref{thm:diagonal-rational}.  For $\chi\in\mathcal C_S$ put $S'=\{i\in S:\chi_i\ne\varepsilon\}$ and $\nu=|S'|$; the condition $\prod\chi_i=\varepsilon$ forces $\nu\ne1$.

\begin{thm}\label{thm:two-diagonal-rational}
Let $X=V(F_1,F_2)\subset\bP^n_\w$ be a nondegenerate pair of diagonal hypersurfaces of the same degree $d$, and $q\equiv1\pmod d$.  Then $\Ztw(\mathcal X,t)$ is a rational function.  Its reciprocal zeros and poles are algebraic numbers of absolute value $q^{w/2}$ with $w\in\ZZ$, $0\le w\le 2(n-2)$; the toric strata contribute the values $q^{j-2}$ with $0\le j\le n$, while the contributions of the character sums have absolute value at most $q^{(n-2)/2}$, with the maximum attained precisely by the full-support tuples with all components nontrivial.
\end{thm}

\begin{proof}
It suffices to write $A_X(q^r)$ as a finite integer combination of geometric sequences $\gamma^r$, since then $\Ztw(\mathcal X,t)$ is a finite product of factors $(1-\gamma t)^{\mp c}$ as in \cref{thm:diagonal-rational}; a rational function with algebraic coefficients and rational Taylor coefficients lies in $\mathbb Q(t)$.  By \cref{prop:hypersurface-count} with $\gcd(k_S,q^r-1)=k_S$ and \cref{lem:pencil-count}, it suffices to treat each term of the lemma over $\F_{q^r}$.

The toric term contributes $\sum_j\binom{|S|-1}{j}(-1)^{|S|-1-j}q^{r(j-2)}$.  For the special-point terms, the Hasse--Davenport relation and $\overline{\chi_i\circ\operatorname{N}}(\gamma)=\overline\chi_i(\gamma)^r$ for $\gamma\in\F_q^\times$ give, exactly as in the proof of \cref{thm:diagonal-rational},
\[
\begin{split}
\Sigma_{S,i_0}(q^r)		&	=(-1)^{|S|-1}\sum_{\chi}\widetilde g_{i_0,\chi}^{\ r},\\
\widetilde g_{i_0,\chi}	&	=\prod_{i\in S\smallsetminus i_0}\bigl(-\overline\chi_i(\gamma_i(P_{i_0}))\,G_q(\chi_i)\bigr),
\end{split}
\]
and the prefactor $(q^r-1)/q^{2r}$ splits each into two geometric terms.

For the middle term, Hasse--Davenport gives
\[
\prod_{i\in S}G_{q^r}(\chi_i\circ\operatorname{N})=(-1)^{|S|}g_{S,\chi}^{\,r}
\]
with $g_{S,\chi}=\prod_{i\in S}(-G_q(\chi_i))$, of absolute value $q^{\nu/2}$.  It remains to stabilise $J_S(\chi\circ\operatorname{N};q^r)$.  Components with $\chi_i=\varepsilon$ contribute the factor $1$ but their points $P_i$ are still excluded, so
\[
J_S(\chi;q^r)=\sum_{\substack{P\in\bP^1(\F_{q^r})\\ \gamma_i(P)\ne0\ \forall i\in S'}}\prod_{i\in S'}\overline{\chi_i\circ\operatorname{N}}\bigl(\gamma_i(P)\bigr)\ -\ \sum_{i_0\in S\smallsetminus S'}\ \prod_{i\in S'}\overline\chi_i\bigl(\gamma_i(P_{i_0})\bigr)^{r},
\]
and the corrections are geometric with terms of absolute value $1$.  For the trivial tuple the main sum is $q^r+1-|S'|$ with $S'=\emptyset$, hence $q^r+1$, and the corrections remove the $|S|$ excluded points.

For $\nu=2$, so $\chi_{i_1}=\overline\chi_{i_2}\ne\varepsilon$, the map
\[
P\mapsto\gamma_{i_2}(P)/\gamma_{i_1}(P)
\]
is a M\"obius isomorphism carrying $\bP^1\smallsetminus\{P_{i_1},P_{i_2}\}$ onto $\Gm$, and the main sum is $\sum_{u\in\F_{q^r}^\times}\chi_{i_1}(u)=0$.  For $\nu\ge3$, all local characters $\overline\chi_i$, $i\in S'$, are nontrivial with trivial product, and by Weil's theory of character sums on curves \cite{Weil1949} there are algebraic numbers $\theta_{S,\chi,1},\dots,\theta_{S,\chi,\nu-2}$, each of absolute value $q^{1/2}$ and each transforming as $\theta\mapsto\theta^r$ under extension, with
\[
	\sum_{\substack{P\in\bP^1(\F_{q^r})\\ \gamma_i(P)\ne0\ \forall i\in S'}}\prod_{i\in S'}\overline{\chi_i\circ\operatorname{N}}\bigl(\gamma_i(P)\bigr)=-\sum_{j=1}^{\nu-2}\theta_{S,\chi,j}^{\ r}.
\]
Assembling, every contribution to $A_X(q^r)$ is an integer multiple of one of $q^{r(j-2)}$, $q^{-r}\widetilde g^{\,r}$, $q^{-2r}\widetilde g^{\,r}$, $q^{-2r}g^{\,r}\zeta^r$ with $\zeta$ of absolute value $1$, or $q^{-2r}g^{\,r}\theta^{\,r}$.  The toric terms have absolute value $q^{j-2}\le q^{n-2}$, attained for $j=n$, $S=T$; among the remaining terms the largest absolute value is $q^{-2}\cdot q^{\nu/2}\cdot q^{1/2}\le q^{(n-2)/2}$, attained only for $\nu=|S|=n+1$.
\end{proof}

\begin{exa}\label{exa:two-conics}
Let $q=13$, $\w=(1,1,1)$, $m=(2,2,2)$, and 
\[
X: x^2+y^2+z^2=y^2-3z^2=0,
\]
 so $a=(1,1,1)$, $b=(0,1,-3)$, with pencil points $P_0=(0:1)$, $P_1=(1:-1)$, $P_2=(3:1)$; the case $b_0=0$ exercises a pencil member missing a variable.  Direct enumeration gives $N_X(T;13)=48$, from $(x,y,z)=(\pm3z,\pm4z,z)$, and $N_X(S;13)=0$ for every proper $S$; thus $A_X(13)=48/12=4$, the four rational intersection points of the two conics.  \cref{lem:pencil-count} reproduces this: with $\rho$ the quadratic character, $G_{13}(\rho)^2=13$, the tuples of $\mathcal C_T$ are the trivial one, with $J=13+1-3=11$, and the three with two components $\rho$, each with $J=-1$ by the evaluation $\sum_s\rho(s^2+bs+c)=-1$ for nonzero discriminant together with the correction at the one excluded point; each special-point sum is $\Sigma=1+13=14$.  Hence 
 \[
 13^2N_X(T;13)=12^3+12\,(-11+3\cdot13)+12^2\cdot42= 169\cdot48.
 \]
\end{exa}

\begin{rem}\label{rem:superelliptic}
The eigenvalues $\theta_{S,\chi,j}$ of \cref{thm:two-diagonal-rational} have a concrete origin: writing $\chi_i=\chi_0^{\mu_i}$ for a character $\chi_0$ of order $N$, they are Frobenius eigenvalues of the $\chi_0$-isotypic part of $H^1$ of the superelliptic curve $y^N=\prod_{i\in S'}(a_is+b_i)^{\mu_i}$, the cyclic cover of the pencil line branched over the special members.  For $\nu=3$ the single eigenvalue is a classical Jacobi sum, and \cref{thm:two-diagonal-rational} degenerates to Gauss--Jacobi data; for $\nu\ge4$ the eigenvalues are genuinely curve-theoretic and no Gauss-sum expression exists in general.  Two limits of the method are worth recording.  Intersections of diagonal hypersurfaces of \emph{unequal} degrees lead to mixed exponential sums in each variable, outside the reach of the elementary analysis; rationality of $\Ztw$ then still holds by general rationality theorems, but without the explicit eigenvalue inventory, and we do not pursue it.  And the recursion of \cref{lem:pencil-count} at the special members reduces the two-equation count to the single-equation count of \cref{lem:diagonal-count}, so the results of this subsection contain those of \cref{subsec:diagonal} as the boundary case.
\end{rem}


\section{Concluding remarks}\label{sec:conclusion}

The guiding theme of this note is that several finite-field counts attached to a weighted projective space, routinely treated as interchangeable, are not.  The coarse point count and the stacky mass coincide and are weight-independent (\cref{cor:coarse-count}, \cref{prop:mass-equals-coarse}); only the orbit count
\[
        A_\w(q)=\sum_{\emptyset\ne S\subseteq T}(q-1)^{|S|-1}\gcd(k_S,q-1)
\]
is sensitive to the weights.  Its significance is not as a counting exercise but as an intrinsic invariant: $A_\w(q)$ is the number of $\F_q$-isomorphism classes of the stack $\cP_\w$, the discrepancy from the coarse count being measured stratum by stratum by the Kummer groups $\F_q^\times/(\F_q^\times)^{k_S}$, and its variation under reduction of weights read off the $\F_q$-cohomology of the $\bbm_d$-gerbe $\cP_\w\to\cP_{\w'}$ (\cref{prop:gerbe-cohomology}, \cref{cor:gerbe-neutral}).

Assembled over extensions, these counts form the twist zeta function $\Ztw(\cP_\w,t)$, whose symmetry theory the paper settles completely.  The function is rational and factors along the isotropy strata, each stratum factor carrying a functional equation centred at its own dimension (\cref{thm:zeta-rational}, \cref{thm:fe}); the sector decomposition of \cref{thm:sector} then computes the multiplicity spectrum in closed form, $b_j=\sum_{n(e)\ge j+1}\varphi(e)$, revealing it to be nonincreasing and independent of $q$.  A single global functional equation therefore exists precisely when the weights share a common prime-to-$p$ part---for reduced $\w$, precisely when every weight is a power of $p$ and the twist theory is trivial (\cref{thm:palindromy-classification}, \cref{cor:no-self-dual}).  The Weil-type symmetry of a weighted projective stack is thus all-or-nothing, in every regime, and the obstruction is the same isotropy that makes $A_\w$ a finer invariant than the coarse count.

The practical upshot is a dictionary.  For questions about rational points on the coarse variety the relevant count is the weight-independent one; for questions about rational representatives, quotient-stack objects, or twisted sectors it is $A_\w(q)$.  The distinction is not academic: for the evaluation codes of \cite{shaska2025} it is the coarse count that governs the parameters, not $A_\w(q)$ (\cref{rem:qwag}), so weighted ambient spaces confer no automatic gain there, whereas $A_\w$ is the correct invariant for twisted objects.

Where the ambient theory closes, the subvariety theory opens.  For diagonal hypersurfaces and same-degree pairs in the split regime, the twist zeta function of the substack is rational with an explicit eigenvalue inventory---Gauss sums, and for intersections the Frobenius eigenvalues of superelliptic curves (\cref{thm:diagonal-rational}, \cref{thm:two-diagonal-rational}); but its spectrum is no longer a $q$-independent combinatorial datum, and whether a functional equation can hold for $\Ztw(\mathcal X,t)$ is a genuine arithmetic question, open even in the diagonal case (\cref{rem:diagonal-general}).  More broadly, $A_\w$ is the local, fixed-field counterpart of the global height counts of \cite{Phillips2024,ESZB} or \cite{shaska-sparsity}; heights, entropy, and the associated zeta functions over finite fields are developed in the companion paper \cite{ss-heights}.


\appendix
\section{Numerical tables}\label{app:tables}

The entries below are obtained directly from the formulas in \ref{subsec:strata-computational}; no search over representatives is required.  They make the numerical part of the paper reproducible and illustrate the size of the smooth and singular contributions over the extensions $q=p^a$.

\begin{table}[hp]
\centering
\caption{Orbit and twist counts for $\w=(2,4,6,10)$ and $q=p^a$, with $p\in\{3,5,7,11\}$ and $1\le a\le5$.}
\label{tab:appendix-24610}
\begin{tabular}{c r r r r}
\toprule
$a$ & $q=p^a$ & $A_{\w}^{\sing}(q)$ & $A_{\w}^{\sm}(q)$ & $A_{\w}(q)$\\
\midrule
1 & 3 & 6 & 74 & 80 \\
2 & 9 & 8 & 1634 & 1642 \\
3 & 27 & 6 & 40874 & 40880 \\
4 & 81 & 16 & 1076162 & 1076178 \\
5 & 243 & 6 & 28816394 & 28816400 \\
\midrule
1 & 5 & 8 & 306 & 314 \\
2 & 25 & 12 & 32546 & 32558 \\
3 & 125 & 8 & 3937746 & 3937754 \\
4 & 625 & 12 & 489063746 & 489063758 \\
5 & 3125 & 8 & 61054693746 & 61054693754 \\
\midrule
1 & 7 & 10 & 794 & 804 \\
2 & 49 & 12 & 240194 & 240206 \\
3 & 343 & 10 & 80943194 & 80943204 \\
4 & 2401 & 20 & 27694108802 & 27694108822 \\
5 & 16807 & 10 & 9495688003994 & 9495688004004 \\
\midrule
1 & 11 & 14 & 2922 & 2936 \\
2 & 121 & 20 & 3572642 & 3572662 \\
3 & 1331 & 14 & 4719441162 & 4719441176 \\
4 & 14641 & 20 & 6277285500482 & 6277285500502 \\
5 & 161051 & 14 & 8354548214002602 & 8354548214002616 \\
\bottomrule
\end{tabular}
\end{table}

\begin{table}[hp]
\centering
\caption{Orbit and twist counts for $\w=(1,2,3,5)$ and $q=p^a$, with $p\in\{3,5,7,11\}$ and $1\le a\le5$.}
\label{tab:appendix-1235}
\begin{tabular}{c r r r r}
\toprule
$a$ & $q=p^a$ & $A_{\w}^{\sing}(q)$ & $A_{\w}^{\sm}(q)$ & $A_{\w}(q)$\\
\midrule
1 & 3 & 4 & 37 & 41 \\
2 & 9 & 4 & 817 & 821 \\
3 & 27 & 4 & 20437 & 20441 \\
4 & 81 & 8 & 538081 & 538089 \\
5 & 243 & 4 & 14408197 & 14408201 \\
\midrule
1 & 5 & 4 & 153 & 157 \\
2 & 25 & 6 & 16273 & 16279 \\
3 & 125 & 4 & 1968873 & 1968877 \\
4 & 625 & 6 & 244531873 & 244531879 \\
5 & 3125 & 4 & 30527346873 & 30527346877 \\
\midrule
1 & 7 & 6 & 397 & 403 \\
2 & 49 & 6 & 120097 & 120103 \\
3 & 343 & 6 & 40471597 & 40471603 \\
4 & 2401 & 10 & 13847054401 & 13847054411 \\
5 & 16807 & 6 & 4747844001997 & 4747844002003 \\
\midrule
1 & 11 & 8 & 1461 & 1469 \\
2 & 121 & 10 & 1786321 & 1786331 \\
3 & 1331 & 8 & 2359720581 & 2359720589 \\
4 & 14641 & 10 & 3138642750241 & 3138642750251 \\
5 & 161051 & 8 & 4177274107001301 & 4177274107001309 \\
\bottomrule
\end{tabular}
\end{table}

\begin{table}[hp]
\centering
\caption{Orbit and twist counts for $\w=(1,6,14,21)$ and $q=p^a$, with $p\in\{3,5,7,11\}$ and $1\le a\le5$.}
\label{tab:appendix-161421}
\begin{tabular}{c r r r r}
\toprule
$a$ & $q=p^a$ & $A_{\w}^{\sing}(q)$ & $A_{\w}^{\sm}(q)$ & $A_{\w}(q)$\\
\midrule
1 & 3 & 13 & 31 & 44 \\
2 & 9 & 37 & 793 & 830 \\
3 & 27 & 109 & 20359 & 20468 \\
4 & 81 & 325 & 537841 & 538166 \\
5 & 243 & 973 & 14407471 & 14408444 \\
\midrule
1 & 5 & 21 & 141 & 162 \\
2 & 25 & 155 & 16201 & 16356 \\
3 & 125 & 501 & 1968501 & 1969002 \\
4 & 625 & 3755 & 244530001 & 244533756 \\
5 & 3125 & 12501 & 30527337501 & 30527350002 \\
\midrule
1 & 7 & 47 & 379 & 426 \\
2 & 49 & 299 & 119953 & 120252 \\
3 & 343 & 2063 & 40470571 & 40472634 \\
4 & 2401 & 14411 & 13847047201 & 13847061612 \\
5 & 16807 & 100847 & 4747843951579 & 4747844052426 \\
\midrule
1 & 11 & 45 & 1431 & 1476 \\
2 & 121 & 731 & 1785961 & 1786692 \\
3 & 1331 & 13323 & 2359716591 & 2359729914 \\
4 & 14641 & 87851 & 3138642706321 & 3138642794172 \\
5 & 161051 & 644205 & 4177274106518151 & 4177274107162356 \\
\bottomrule
\end{tabular}
\end{table}

\end{document}